\documentclass[11pt]{article}
\def\be{\begin{equation}}
\def\ee{\end{equation}}
\def\beq{\arraycolsep=1.5pt\begin{eqnarray}}
\def\eeq{\end{eqnarray}}
\def\R{I\!\!R}
\def\D{\tilde{D}}
\def\t{\tilde}
\def\div{{\,\rm div \,}}
\def\IS{{\,\rm IS \,}}
\def\kf{{\,\rm kf \,}}
\def\KF{{\,\rm KF \,}}
\def\sym{{\,\rm sym \,}}
\def\is{{\,\rm is \,}}
\def\dist{{\,\rm dist \,}}
\def\SO{{\,\rm SO \,}}

\def\ib{{\,\rm ib \,}}
\def\Ric{{\,\rm Ric \,}}
\def\V{{\,\cal V \,}}
\def\ctg{{\,\rm ctg \,}}
\def\Vol{{\,\rm Vol \,}}

\def\dim{{\,\rm dim \,}}
\def\Kill{{\,\rm KF \,}}
\def\S{{\,\rm S \,}}
\def\T(M){{\,\rm T(M) \,}}
\def\Om{\Omega}
\def\qfq{{\quad\mbox{for}\quad}}
\def\Ga{\Gamma}
\def\<{\left<}
\def\>{\right>}
\def\si{\sigma}
\def\lam{\lambda}
\def\A{{\cal A}}

\def\pl{{\partial}}
\def\F{{\cal F}}
\def\i{{\,\rm i \,}}
\def\B{{\,\cal B \,}}
\def\tr{{\,\rm tr \,}}
\def\a{{\alpha}}
\def\b{{\beta}}
\def\X{{\cal X}}
\def\Q{{\cal Q}}

\def\var{\varphi}
\usepackage{amsfonts}
\usepackage{latexsym, amssymb, amsmath, amscd, amsfonts, epsfig, graphicx}
\usepackage{mathrsfs}

\newfont{\Blackboard}{msbm10 scaled 1200}

\newfont{\roma}{cmr10 scaled 1200}

\newtheorem{thm}{{}\hskip\parindent Theorem}[section]
\newtheorem{lem}{{}\hskip\parindent Lemma}[section]
\newtheorem{pro}{{}\hskip\parindent Proposition}[section]
\newtheorem{exl}{{}\hskip\parindent Example}[section]
\newtheorem{cor}{{}\hskip\parindent Corollary}[section]

\newtheorem{rem}{{}\hskip\parindent Remark}[section]

\def\beq{\arraycolsep=1.5pt\begin{eqnarray}}
\def\eeq{\end{eqnarray}}

\renewcommand{\theequation}{\arabic{\section}.\arabic{\equation}}
\numberwithin{equation}{section}

\large
\title{Space of Infinitesimal Isometries and Bending of Shells}
\date{}
\author{Peng-Fei Yao\\[0.3cm]
Key Laboratory of Systems and Control \\
Institute of Systems Science,
Academy of Mathematics and Systems Science\\
Chinese Academy of Sciences, Beijing 100080, P. R. China\\
e-mail: pfyao@iss.ac.cn \\[0.3cm]}

\begin{document}
\maketitle \footnote{This work is  supported by the National Science
Foundation of China, grants  no. 60821091 and no. 61174083. }
\begin{quote}
\begin{small}
{\bf Abstract} \,\,\,We discuss infinitesimal isometries of the
middle  surfaces and present  some characteristic conditions for a
function to be the normal component of an infinitesimal isometry.
Our results show that those characteristic conditions depend on
the Gaussian curvature of the middle surfaces: Normal components
of infinitesimal isometries satisfy an elliptic problem, or a
parabolic one, or a hyperbolic one according to the middle surface
being elliptic, or parabolic, or hyperbolic, respectively. In
those cases, a problem of determining an infinitesimal isometry is
changed into that of $1$-dimension. Then we apply those results to
the energy functionals of bending of shells which has been
obtained as two-dimensional problems by the limit theory of
$\Ga$-convergence from the three-dimensional nonlinear elasticity.
Therefore the limit theory of $\Ga$-convergence reduces to be a
one-dimensional problem in the those cases.
\\[3mm]
\end{small}
\end{quote}

\tableofcontents

\setcounter{equation}{0}
\section{Introduction}
\def\theequation{1.\arabic{equation}}
\hskip\parindent Let $M\subset\R^3$ be a smooth surface and let
$\Om\subset M$ be a bounded, open set. A map $V:$
$\Om\rightarrow\R^3$ is said to be an infinitesimal isometry on
$\Om$ if $$\<\hat D_XV,X\>=0\qfq X\in M_x,\,\,x\in\Om,$$ where
$\<\cdot,\cdot\>$ is the Euclidean metric of $\R^3$ and $\hat D$
denotes the covariant differential of the Euclidean space $\R^3.$
We denote by $\IS^1(\Om,\R^3)$ all $H^1$ infinitesimal isometries
on $\Om.$

The study of infinitesimal isometries has been a long history, see
\cite{Al, CoVo, Pog1, Pog, Sp} and many others. Their purposes
were to establish "infinitesimal rigidity" for some closed
surfaces and their interests were not on the structure of
infinitesimal isometries themselves. For a detail survey along
this direction, we refer to \cite{Sp}.

Our interests in the space $\IS(\Om,\R^3)$ of infinitesimal
isometries are motivated by the recent lower dimensional models for
thin structures(such as membrane and shells) through
$\Ga$-convergence. This approach has lead to the derivation of a
hierarchy of limiting theories and provides a rigorous justification
of convergence of three-dimensional minimizers to minimizers of
suitable lower dimensional limit energies.

Given a 2-dimensional surface $\Om$, consider a shell $S^h$ of
middle surface $\Om$ and thickness $h$, and associate to its
deformation $u$ the scaled per unit thickness three dimensional
nonlinear elastic energy $E(u,S^h).$  The $\Ga$-limit $I^\b$ of the
energies
$$h^{-\b}E(u,S^h)$$ are identified as $h\rightarrow0$ for a given
scaling $\b\geq0.$

When $\Om$ is a subset of $\R^2$ (i.e., a plate), such
$\Ga$-convergence was first established by \cite{LeRa} for $\b=0,$
then by \cite{FrMoMu1, FrMoMu2} for $\b\geq2$ (also see \cite{PaO}
for the results for $\b=2$ under additional conditions). In case
of $0<\b<5/3$ the convergence was obtained by \cite{CoMa}, see
also \cite{Conti}. Other significant results for plates concern
the derivation of limit theories for incompressible materials
\cite{CoDo, CoDo1, Trabe}, for heterogeneous materials
\cite{Schmidt}, and through establishing convergence of
equilibria, rather than strict minimizers \cite{Monn, MuPa}.

When $\Om$ is an any surface, the first result by \cite{LeRa1}
relates to scaling $\b=0$ and models membrane shells. The limit
energy $I^0$ depends  only on the stretching and shearing produced
by the deformation of the middle surface $\Om.$ Then the limit
energy $I^2$ in the case of $\b=2$ was given by \cite{FrMoMu}. This
scaling corresponds to a flexural shell model, where the only
admissible deformations are those preserving the metric of $\Om.$
Then the energy $I^2$ depends on the change of curvature by the
deformation. The limit energies $I^\b$ are obtained by \cite{LeMoPa}
for scaling $\b\geq4.$ Based on some quantitative rigidity estimate
due to \cite{FrMoMu1},  \cite{LeMoPa} demonstrates that the first
term in the expansion $u-R$, in terms of $h$, belongs to the space
$\IS(\Om,\R^3)$ of infinitesimal isometries. That means that there
is no first order change in the induced metric of the meddle surface
$\Om$. The corresponding limit energy $I^\b$ consists of the bending
energy which is given by the first order change of the second
fundamental form of $\Om$ for $\b>4.$  In the case of $\b=4$,
\cite{LeMoPa} also shows that, if the middle surface $\Om$ is
approximately robust, the $\Ga$-limit is still a bending term.
Moreover, in the scaling regime of $2<\b<4$  the limit $I^\b$ is
given  by \cite{LeMoPa2} which reduces to be the pure bending energy
again.

As shown by \cite{LeMoPa2, LeMoPa}, the limit energy functionals
$I^\b$ of the $\Ga$-convergence for all $\b>2$ are over the space
$\IS(\Om,\R^3)$ of infinitesimal isometries. Then  the space
$\IS(\Om,\R^3)$ naturally plays a crucial role in the analysis of
shells. The aim of the present paper is to understand the space
$\IS(\Om,\R^3).$

We now give heuristic overview of our results, whose precise
formulations will be presented in the sections later. Let $N$ be the
unit normal field of surface $M$ and  let $\X(\Om)$ be all vector
fields on $\Om$. For $V\in H^1(\Om,\R^3)$, we decompose as
$$V=W+wN\qfq W\in\X(\Om),\quad w\in H^1(\Om).$$ We look for
conditions on  functions $w$ such that there are  vector fields
$W\in\X(\Om)$ to guarantee  $V\in \IS^1(\Om,\R^3).$

First, Section 2 is devoted to treating the structure of
$V=W\in\IS^1(\Om,\R^3)$ corresponding to the zero normal component
$w=0.$ Such an infinitesimal isometry is said to be a Killing field.
Through there are rich results on Killing fields (\cite{Pe}), we
focus on the relations between a Killing field and the Gaussian
curvature function. In particular, we show that the dimension of the
Killing field space is $3$ if $\Om$ is of constant curvature and is
not larger than $1$ in the case of non-constant curvature
(Corollaries \ref{nc2.1}, \ref{c2.2}, and Theorem \ref{nt2.3}).
Furthermore the explicit formulas of Killing fields are given in
terms of the Gaussian curvature function (Theorem \ref{t2.2}).

Let $H^1_\is(\Om)$ denote all functions $w\in H^1(\Om)$ such that
there are  $W\in\X(\Om)$, which are perpendicular to all Killing
fields in $H^1(\Om,\R^3)$, to ensure that
$V=W+wN\in\IS^1(\Om,\R^3).$  Section 3 shows that $w\in
H^1_\is(\Om)$ if and only if $w$ satisfies an equation
(\ref{3.12f2}) (Theorem \ref{th3.1}).

The type of the equation (\ref{3.12f2}) is subject to the Gaussian
curvature function: It is elliptic, or parabolic, or hyperbolic
according to ellipticity, or parabolicity, or hyperbolicity of the
middle surface $\Om$, respectively. The three cases are studied,
respectively, in Sections 4, 5, and 6.  Our results show that the
problem to determine whether $w\in H^1_\is(\Om)$ is actually that
of $1$-dimension in the above three types, respectively.

As a consequence of those theories, we present a condition for the
middle surface $\Om$ which can guarantee that $H^1_\is(\Om)\cap
C^\infty(\Om)$ is dense in $H^1_\is(\Om)$ in the norm of $H^1(\Om)$
(Theorems \ref{tn4.3} and \ref{nt6.3}): There is a point $o\in\Om$
such that $\Om$ is star-shaped with respect to $o$ and
$$\Om\subset\exp_o\Sigma(o),$$ where $\exp_o\Sigma(o)$ is the interior of the cut locus of $o.$
Such an issue is actually not trivial. In general, even though $\Om$
is elliptic, an element $V\in \IS^1(\Om,\R^3)$ may not be
approximated by smooth infinitesimal isometries. An interesting
example, discovered by \cite{CoVo} (also see \cite{Sp}), is a closed
smooth surface of non-negative curvature for which the infinitesimal
rigidity holds true: All $C^\infty$ infinitesimal isometries are
trivial.  But there is a $C^2$ non-trivial infinitesimal isometry.
Therefore $H^1_\is(\Om)\cap C^\infty(\Om)$ is not dense in
$H^1_\is(\Om)$ for this surface.

In Section 7 we apply the above theories to the limit energy $I^\b$
of $\Ga$-convergence for the scaling $\b>2.$ Then the limit energy
functional is changed into a one-dimensional formula over a function
space with one variable (Theorem \ref{t8.1}). In particular, we
present the explicit formulas of the limit energy functionals for a
spherical shell (Theorem \ref{t8.2}) and a cylinder shell (Theorem
\ref{t8.3}), respectively, under the nonlinear isotropic materials.

Here we do not use the traditional methods, adopted in the classical
linear thin shell theories. Their starting point is to assume that
the middle surface is given by a coordinate path: $\Om$ is the image
in $\R^3$ of a smooth map defined on a connected domain of $\R^2$,
rooted from classical differential geometry. The classical models
use the traditional geometry and end up with highly complicated
resultant equations. In these, the explicit presence of the
Christoffel symbols, makes some necessary computations too
complicated. We view the middle surface $\Om$ as a $2$-dimensional
Riemannian manifold with the induced metric to make everything
coordinates free as far as possible. When necessary, some special
coordinates are chosen to simplify computations as in modelling and
control for the classical thin shells, see \cite{Chai, Chai1, ChGu,
ChGuYa, ChLi, ChLi1, FeFe, LaTr, LiYa,  Yao2, Yao, Yao3, Yao4} and
many others.

\section{Killing Fields in Dimension $2$}
\def\theequation{2.\arabic{equation}}
\hskip\parindent We shall present explicit formulas of Killing
fields in terms of the Gaussian curvature function (Proposition
\ref{p2.1} and Theorem \ref{t2.2}).

 Let $M\subset\R^3$ be a smooth surface with the induced
metric $g$ from the Euclidean metric of $\R^3.$ Let $\Om\subset M$
be an open set with smooth boundary $\Ga.$ Denote by $\X(\Om)$ all
vector fields on $\Om$. Let $W\in\X(\Om)$ be given. Let $\a(t)$ be
the 1-parameter group generated by $W$, i.e.,
\be\dot{\a}(t,x)=W(\a(t,x)),\quad \a(0,x)=x\qfq x\in
\Om.\label{n2.1}\ee $W$ is said to be a Killing field on $\Om$ if
and only if $\a(t)$ are local isometries.  It is easy to check that
$W$ is a Killing field on $\Om$ if and only if $V=W$ is a $C^\infty$
infinitesimal isometry on $\Om$, that is, \be DW(X,Y)+DW(Y,X)=0\qfq
X,Y\in M_x,\quad\mbox{for all}\quad x\in\Om,\label{2.1}\ee where $D$
is the Levi-Civita connection of the induced metric $g$ and  $DW$ is
the covariant differential of vector field $W$.

Let \be\Kill(\Om,T)=\{\,\mbox{all $C^\infty$ Killing fields on
$\Om$}\,\}.\label{2.10*}\ee Then (\cite{Pe})$$\dim\Kill(M,T)\leq
3.$$

For the purpose of application to bending of shells here, we only
need to consider $H^1$ Killing fields. To this end, we introduce
some common  notions in Riemannian geometry. Let $T^k(\Om)$ be all
$k$-th order tensor fields on $\Om$ where $k$ is a nonnegative
integer. In particular, $T^0(\Om)$ is all functions on $\Om$ and
$T(\Om)=\X(\Om).$ For each $x\in M$,
 the k-th order tensor space $T^k_x$ on $M_x$ is an inner product space
defined as follows. Let $e_1$,  $e_2$ be an orthonormal basis of
$M_x$. For any $\a$, $\b\in T^k_x$, $x\in M$, the inner product is
given by \be
 \<\a,\,\b\>_{T^k_x}=\sum_{i_1=1,\cdots,i_k=1}^2\a(e_{i_1},\cdots,e_{i_k})
\b(e_{i_1},\cdots,e_{i_k})\quad\mbox{at}\quad x. \label{8.4}\ee Note
that the right hand side of (\ref{8.4}) is free of choice of
orthonormal bases. In particular, for $k=1$ the definition
(\ref{8.4}) becomes
$$ g(\a,\b)=\<\a,\b\>_{T_x^1}=\<\a,\b\>, \quad\forall\,\,\a,\,\b\in M_x,
$$ that is, the induced inner product of $M_x$ of $M$ from $\R^3$.  Let $L^2(\Om,T^k)$
be the Soblev spaces of all $k$-th order tensor fields on $\Om$ with
inner products
$$(T_1, T_2)_{L^2(\Om,T^k)}=\sum_{i=0}^k\int_\Om\<T_1,T_2\>_{T^k_x}dg\qfq T_1,\,\,T_2\in L^2(\Om,T^k).$$
Let $$H^1(\Om,T)=\{\,\,W\,\,|\,\,W\in\X(\Om),\,\,W\in
L^2(\Om,T),\,\,DW\in L^2(\Om,T^2)\,\,\}$$ with norm
$$\|W\|_{H^1(\Om,T)}=\Big(\|w\|^2_{L^2(\Om,T)}+\|DW\|^2_{L^2(\Om,T^2)}\Big)^{1/2}.$$

Denote  by $H^1_\kf(\Om,T)$ all $W\in H^1(\Om,T)$ with the
relations (\ref{2.1}) being true for almost everywhere on $\Om$
and  the norm of $H^1(\Om,T).$ Then
\be\KF(M,T)\subset\KF(\Om,T)\subset H^1_\kf(\Om,T).\label{n2.5}\ee

We have
\begin{thm}\label{tn3.1} Let $\Om\subset M$ be an open set. Then
$$H^1_\kf(\Om,T)=\KF(\Om,T)$$ and
\be\dim\KF(M,T)\leq\dim H^1_\kf(\Om,T)\leq3.\label{n2.6}\ee
\end{thm}

{\bf Proof}\,\,\,Let $W\in H^1_\kf(\Om,T).$ By Lemma 4.4 in
 \cite{Yao}, we have
 \be {\bf \Delta}W=2\kappa W\qfq x\in \Om,\label{n2.77}\ee
 where ${\bf \Delta}$ is the Hodge-Laplace operator in the metric $g$
  and $\kappa$ is the Gaussian curvature function of $M$.  Then the
  ellipticity of the operator ${\bf\Delta}$ implies that $W$ is
  $C^\infty$ on $\Om$, that is, $W\in\KF(\Om,T).$ Then the left hand side of the inequality
  (\ref{n2.6}) follows from (\ref{n2.5}). Moreover, the right hand side of the inequality (\ref{n2.6})
 is given by the equations (\ref{2.5})-(\ref{2.8*}) later. \hfill$\Box$

\begin{rem}In general, $\dim\KF(M,T)\not=\dim H^1_\kf(\Om,T),$ see Corollary $\ref{c2.2}$ and Example
$\ref{e3.1}$
later.
\end{rem}

 Let $o\in M$ be fixed and let $\exp_o:$ $M_o\rightarrow M$
be the exponential map in the metric $g.$ For any $v\in M_o$ with
$|v|=1,$ then there is a unique $t_0(v)>0$ (or $t_0(v)=\infty$) such
that the normal geodesic $\gamma(t)=\exp_otv$ is the shortest on the
interval $[0,t_0].$  Let
$$C(o)=\{\,\,t_0(v)v\,\,|\,\,v\in M_o,\,\,|v|=1\,\,\},\quad
\Sigma(o)=\{\,\,tv\,\,|\,\,v\in M_o,\,\,|v|=1,\,\,0\leq t<
t_0(v)\,\,\}.$$ The set $\exp_oC(o)\subset M$ is said to be the
cut locus of $o$ and the set $\exp_o\Sigma(o)\subset M$ is called
the interior of the cut locus of $o.$ Then
$$M=\exp_o\Sigma(o)\cap\exp_oC(o).$$ Furthermore, $\exp_o:$ $\Sigma(o)\rightarrow
\exp_o\Sigma(o)$ is a diffeomorphism and $C(o)$ is a zero measure
set on $M_o$. Then $\exp_o C(o)$ is a zero measure set on $M$
since it is the image of the zero measure set $C(o)$, that is,
$\exp_o\Sigma(o)$ is $M$ minus a zero measure set.

We introduce the polar coordinate system at $o\in M$ as follows.
Let $e_1,$ $e_2$ be an orthonormal basis of $M_o$. Set
\be\si(\theta)=\cos\theta e_1+\sin\theta e_2\qfq
\theta\in[0,2\pi).\label{2.2**}\ee Consider a family of two
parameter curves on $M$ given by
$$ {\cal F}(t,\theta)=\exp_ot\si(\theta)\qfq t\si(\theta)\in\Sigma(o).$$
Then \be \pl t=\frac{\pl}{\pl t}{\cal
F}(t,\theta)=\exp_{o*}\si(\theta),\quad
\pl\theta=\frac{\pl}{\pl\theta}\F(\,t\,,\theta)=\,t\,
\exp_{o*}\dot{\si}(\theta).\label{2.4**}\ee In particular,
$$g=d\,t\,^2+f^2(\,t\,,\theta)d\theta^2\qfq x=\exp_o\,t\,\si(\theta)\in\exp_o\Sigma(o),$$
where $f(t,\theta)$ is the solution to the problem
 \be\left\{\begin{array}{l}f_{tt}(t,\theta)+\kappa(t,\theta)
f(t,\theta)=0,\\ f(0,\theta)=0,\quad
f_t(0,\theta)=1,\end{array}\right.\label{2.2}\ee where $\kappa$ is
the Gaussian curvature function on $M$ and
$\kappa(t,\theta)=\kappa(\F(t,\theta)).$

 Let \be T=\pl t,\quad E=\frac{1}{f}\pl\theta\qfq
 x\in\exp_o\Sigma(o)-\{o\}.\label{n2.11}\ee
 Then $T,$ $E$ is a frame field on $\exp_o\Sigma(o)-\{o\}.$
 Let $D$ denote the Livi-Civita connection of the induced metric $g$ on $M.$ We have \be
D_TT=0,\quad D_TE=0,\quad D_ET=\frac{f_t}{f}E,\quad
D_EE=-\frac{f_t}{f}T\ee for $x\in\exp_o\Sigma(o)-\{o\}.$

 Let
$$W=\var(t,\theta) T+\phi(t,\theta) E,$$ where $\var=\<W,T\>$ and
$\phi=\<W, E\>.$ Then \be D_TW=\var_tT+\phi_tE,\quad
D_EW=\frac{1}{f}(\var_\theta-f_t\phi)T+\frac{1}{f}(\phi_\theta+f_t\var)E.\label{n2.7}\ee
Then the relation
$(\ref{2.1})$ is equivalent to \be\left\{\begin{array}{l} \var_t(t,\theta)=0,\\
f\phi_t-f_t\phi+\var_\theta=0,\\
\phi_\theta+f_t\var=0.
\end{array}\right.\label{2.5}\ee

The first equation in (\ref{2.5}) yields
\be\var(t,\theta)=\<W,\si(\theta)\>(o).\label{solu1}\ee Next, we
calculate the first order derivative of the second equation in
(\ref{2.5}) with respect to $t$ and use the equations
(\ref{solu1}) and (\ref{2.2}) to have \be
\left\{\begin{array}{l}\phi_{tt}+\kappa\phi=0,\\
\phi(0)=\<W,\dot{\si}(\theta)\>,\quad\phi_t(0)=DW(e_2,e_1).\end{array}\right.\label{2.8*}\ee

It follows from (\ref{solu1}) and (\ref{2.8*}) that

\begin{pro}\label{p2.1} Let $(M,g)$ be of constant curvature $\kappa.$  Let $\Om\subset M$ be
an open set and let $o\in\Om.$ For $W\in\KF(\Om,T)$, we have
 \be
 W=\<W(o),\si(\theta)\>T+(af+f_t\<W(o),\dot{\si}(\theta)\>)E\label{2.4}\ee
for $x$ in a neighborhood of $o$,  where $f$ is given by
$(\ref{2.2})$ and $a=DW(e_2,e_1).$
\end{pro}

Let $M$ be a sphere with curvature $\kappa>0.$ Let $W\in\KF(M,T).$
Then $f=\dfrac{1}{\sqrt{\kappa}}\sin\sqrt{\kappa}t$ for
$t\in(0,\dfrac{\pi}{\sqrt{\pi}})$ and the formula (\ref{2.4})
holds true for all $x\in\exp_o\Sigma(o).$ It is easy to check that
$W$, given by (\ref{2.4}), is $C^\infty$ at $\exp_oC(o)$, which is
the antipodal point of $o.$  Then $$\dim\KF(M,T)=3.$$

It follows Theorem \ref{tn3.1} that

\begin{cor}\label{nc2.1}  Let $M\subset\R^3$ be a closed sphere of constant curvature $\kappa>0$ and let
 $\Om\subset M$ be an open set. Then
$$\dim H^1_\kf(\Om,T)=3.$$
\end{cor}

Let $o\in\Om$ be given. $\Om$ is said to be star-shaped with respect
to $o$ if for any $x\in\Om$ there is a shortest geodesic contained
in $\Om$ connecting $x$ and $o.$ Since for any $W(o)\in M_o$ and a
number $a=DW(e_2,e_1)$ given, the problems (\ref{solu1}) and
(\ref{2.8*}) have solutions for all $x\in\exp_o\Sigma(o),$ it follows
that

\begin{cor} \label{c2.2}Let $M\subset\R^3$ be a surface with zero curvature. Let $\Om\subset M$ be star-shaped with respect to
$o\in\Om$ and \be\Om\subset\exp_o\Sigma(o).\label{new2.13}\ee Then
\be\dim H^1_\kf(\Om,T)=3.\label{m2.19}\ee
\end{cor}

\begin{rem} The condition $(\ref{new2.13})$ is necessary for the equation $(\ref{m2.19})$. This is
because a vector field $W$, which is given by  $(\ref{2.4})$, can
not guarantee  $W$ is $C^\infty$ on $\Om\cap\exp_oC(o),$ see
Example $\ref{e3.1}$ below.
\end{rem}

\begin{exl}\label{e3.1} Consider a cylinder
$$M=\{\,\,(x,z)\,\,|\,\,x=(x_1,x_2)\in\R^2,\,\,|x|=1,\,\,z\in\R\,\,\}.$$
Let $b>0$ and let $$\Om=\{\,\,(x,z)\,\,|\,\,|x|=1,\,\,|z|<b\,\,\}.$$
Let $o=(1,0,0).$ Then
$$\Om\cap\exp_oC(o)=\{\,\,(-1,0,z)\,\,|\,\,|z|<b\,\,\},\quad f(t,\theta)=t,$$ and the
vector field
$$W=tE$$ is not well defined on $\Om\cap\exp_oC(o).$ In this case,
it is easy to check that
$$\dim H^1_\kf(\Om,T)=2.$$   In particular, $$\dim\KF(M,T)=2.$$
\end{exl}

\begin{lem} \label{l2}Let $\kappa$ be the Gaussian curvature function of $M$ and let $W$ be a Killing field on $\Om$. Then
\be \<\nabla\kappa,W\>=0\qfq x\in \Om,\label{2.2*}\ee \be
D^2\kappa(\nabla\kappa,W)=0\qfq x\in\Om.\label{2.3*}\ee
\end{lem}

{\bf Proof}\,\,\,Let $o\in \Om$ be given. We have \beq\kappa_{\theta
t}&&=(f\<\nabla\kappa,E\>)_t=f_t\<\nabla\kappa,E\>+fD^2\kappa(E,T),\nonumber\eeq
which gives
$$\kappa_{\theta t}(0)=\<\nabla\kappa,\dot{\si}(\theta)\>.$$

Let $\var$ and $\phi$ be given by $(\ref{solu1})$ and
$(\ref{2.8*}),$ respectively. Using the equations (\ref{2.2}) and
(\ref{2.8*}) and the third equation in (\ref{2.5}), we obtain \beq
0&&=\phi_\theta^{(3)}(0)+f^{(4)}(0)\var(0)=-\kappa_{\theta
t}(0)\phi(0)-\kappa_t(0)\phi_\theta(0)-\kappa(0)\phi_{\theta t}(0)-2\kappa_t(0)\var(0)\nonumber\\
&&=-\<\nabla\kappa,W\>(o),\nonumber\eeq
 that is, the formula (\ref{2.2*}) is true at $o$, where the
following formula is used
$$\phi_\theta(0)=\<W,\ddot{\si}(\theta)\>=-\<W,\si(\theta)\>,$$
$$\phi_{\theta t}(0)=[DW(\dot{\si}(\theta),\si(\theta))]_\theta=-DW(\si(\theta),\si(\theta))+
DW(\dot{\si}(\theta),\dot{\si}(\theta))=0.$$ Since $o\in \Om$ can be
any point, the formula (\ref{2.2*}) follows.

Finally using (\ref{2.2*}), we have
$$0=\nabla\kappa\<\nabla\kappa,W\>=D^2\kappa(\nabla\kappa,W)+
\<\nabla\kappa,D_{\nabla\kappa}W\>=D^2\kappa(\nabla\kappa,W).\quad\quad\quad\quad\quad\quad\hfill\Box$$

Let $(M,g)$ be orientable. Let $X$ be a vector field on $\Om$. We
define a vector field $QX$ on $\Om$ by\be
QX=\<X,e_2\>e_1-\<X,e_1\>e_2\qfq x\in \Om,\label{2.9*}\ee where
$e_1,$ $e_2$ is an orthonormal basis of $M_x$ with an positive
orientation. It is easy to check that the vector field $QX$ is well
defined.

 We
have

\begin{thm}\label{t2.2} Let $(M,g)$ be orientated and let
$\kappa$ be the Gaussian curvature function.                   Let
$\Om\subset M$ be a connected open set and let
$$|\nabla\kappa|>0\qfq x\in\Om.$$ Then \be\dim H^1_\kf(\Om,T)=1\label{n2.21}\ee
 holds true if and only if the following formulas are true \be
D^2\kappa(\nabla\kappa,Q\nabla\kappa)=0\qfq x\in\Om,\label{2.12*}\ee
\be \<Q\nabla\kappa,\nabla\Delta\kappa\>=0\qfq
x\in\Om,\label{2.13*}\ee where $\Delta$ is the Laplacian of the
metric $g$. Moreover, $W\in H^1_\kf(\Om,T)$ has a formula \be
W=ce^{h_0}Q\nabla\kappa,\ee where $c$ is a constant and $h_0$ is a
solution to the problem \be\nabla
h=\frac{|\nabla\kappa|^2\Delta\kappa-2D^2\kappa(\nabla\kappa,\nabla\kappa)}{|\nabla\kappa|^4}\nabla\kappa.\label{hn}\ee
\end{thm}

{\bf Proof}\,\,\, By Lemma \ref{l2}, we look for Killing fields in
the form \be W=e^hQ\nabla\kappa,\label{2.14*}\ee where $h$ is a
function on $\Om$.

Let $o\in\Om$ be given. Let $e_1,$ $e_2$ be an orthonormal basis of
$M_o$ with an positive orientation. Let $\si(\theta)$ be given by
(\ref{2.2**}). Then
$$T=T(t,\theta),\quad E=E(t,\theta)$$ forms an orthonormal basis of $M_{\F(t,\theta)}$ with the positive
orientation for $\F(t,\theta)\in\Om$. For convenience, we denote
$$E_1=T,\quad E_2=E,$$ and
$$p_i=\<\nabla p,E_i\>,\quad p_{ij}=D^2p(E_i,E_j),\quad
p_{ijk}=D^3p(E_i,E_j,E_k),$$ for $i,$ $j,$ $k=1,$ $2$, where $p$ is
a function on $\Om$. Then
$$W=e^h\kappa_2E_1-e^h\kappa_1E_2.$$ Using the relations
$$D_{E_1}E_1=D_{E_1}E_2=0,\quad D_{E_2}E_1=\frac{f_t}{f}E_2,\quad
D_{E_2}E_2=-\frac{f_t}{f}E_1,$$ we obtain
$$\left\{\begin{array}{l} e^{-h}D_{E_1}W=(h_1\kappa_2+\kappa_{12})E_1-(h_1\kappa_1+\kappa_{11})E_2,\\
e^{-h}D_{E_2}W=(h_2\kappa_2+\kappa_{22})E_1-(h_2\kappa_1+\kappa_{12})E_2.
\end{array}\right.$$  It follows from the above formulas that $W$ is
a Killing field if and only if
\be\left\{\begin{array}{l}h_1\kappa_2+\kappa_{12}=0,\\
h_2\kappa_1+\kappa_{12}=0,\\
-h_1\kappa_1+h_2\kappa_2+\kappa_{22}-\kappa_{11}=0.\end{array}\right.\label{2.15*}\ee
Then the formula (\ref{2.14*}) defines a Killing field if and only
if there is a solution $h$ to the problem (\ref{2.15*}).

Now we solve the problem
\be\left\{\begin{array}{l}h_1\kappa_2+h_2\kappa_1+2\kappa_{12}=0,\\
-h_1\kappa_1+h_2\kappa_2+\kappa_{22}-\kappa_{11}=0,\end{array}\right.\label{2.16**}\ee
to
have\be\left(\begin{array}{c}h_1\\h_2\end{array}\right)=-\frac{1}{|\nabla\kappa|^2}
\left(\begin{array}{c}2\kappa_2\kappa_{12}+\kappa_1(\kappa_{11}-\kappa_{22})\\
2\kappa_1\kappa_{12}-\kappa_2(\kappa_{11}-\kappa_{22})\end{array}\right).\label{2.16*}\ee
It is easy to check that a solution of (\ref{2.15*}) is a solution
of (\ref{2.16*}) if and only if the formula (\ref{2.12*}) holds.
Then the formula (\ref{2.14*}) defines a Killing field if and only
if the formula (\ref{2.12*}) holds true and  $h$ satisfies the
problem (\ref{2.16*}).

Next, let us show that there is a solution to the problem
(\ref{2.16*}) if and only if the formula (\ref{2.13*}) holds true.
To this end, we let
$$X=h_1E_1+h_2E_2,$$ where $(h_1,h_2)^\tau$ is given by the formula
(\ref{2.16*}). We review the vector field $X$ as a $1$-form. It
follows from the Poincare lemma that the problem (\ref{2.16*}) has a
solution if and only if
$$dX=0\qfq x\in\Om,$$ since $\Om$ is star-shaped with respect to
$o$, where $d$ denotes the exterior differentiation.

We assume that the formula (\ref{2.12*}) holds true to compute $dX$.
By \cite{Wu}, we have \beq dX&&=E_1\wedge D_{E_1}X+E_2\wedge
D_{E_2}X=[E_1(h_2)-E_2(h_1)+h_2f_t/f]E_1\wedge E_2.\nonumber\eeq
Using (\ref{2.15*}) and (\ref{2.16*}), we obtain \be
E_2(h_1)\kappa_2^2=\kappa_{12}\kappa_{22}-\kappa_2\kappa_{122}+[h_2|\nabla\kappa|^2+\kappa_1
\kappa_{12}]f_t/f, \label{2.17*}\ee \beq
E_2(h_1)(\kappa_2^2-\kappa_1^2)&&=2h_1\kappa_1\kappa_{12}+\kappa_{12}(\kappa_{11}+\kappa_{22})
+\kappa_1(\kappa_{112}-\kappa_{222})\nonumber\\
&&\quad+\kappa_2(\kappa_{11}-\kappa_{22})f_t/f.\label{2.18**}\eeq It
follows from (\ref{2.17*}), (\ref{2.18**}), and (\ref{2.16*}) that
\beq
E_2(h_1)|\nabla\kappa|^2&&=-2h_1\kappa_1\kappa_{12}-\kappa_{12}(\kappa_{11}-\kappa_{22})-2\kappa_2\kappa_{122}\nonumber\\
&&\quad-\kappa_1(\kappa_{112}-\kappa_{222})+h_2|\nabla\kappa|^2f_t/f.\label{2.19*}\eeq
Similarly, we have \be
E_1(h_2)|\nabla\kappa|^2=-2h_2\kappa_2\kappa_{12}+\kappa_{12}(\kappa_{11}-\kappa_{22})-2\kappa_1\kappa_{121}
+\kappa_2(\kappa_{111}-\kappa_{221}).\label{2.20*}\ee

From (\ref{2.19*}), (\ref{2.20*}), and (\ref{2.16**}), we obtain
\beq
&&[E_1(h_2)-E_2(h_1)+h_2f_t/f]|\nabla\kappa|^2\nonumber\\
&&=2\kappa_{12}[h_1\kappa_1-h_2\kappa_2+\kappa_{11}-\kappa_{22}]-\kappa_1(\kappa_{112}+\kappa_{222})+\kappa_2(\kappa_{111}+\kappa_{221})\nonumber\\
&&=\<Q\nabla\kappa,\nabla\Delta \kappa\>,\label{2.22}\eeq where the
following formulas have been used
$$\kappa_{121}=\kappa_{112}+\kappa_2\kappa,\quad
\kappa_{122}=\kappa_{212}=\kappa_{221}+\kappa_1\kappa,$$
$$E_2(\Delta\kappa)=\kappa_{112}+\kappa_{222}+2(\kappa_{12}-\kappa_{12})f_t/f=\kappa_{112}+\kappa_{222},\quad
E_1(\Delta\kappa)=\kappa_{111}+\kappa_{221}.$$

To complete the proof, it remains to show that  the formula
(\ref{2.16*}) is the same as (\ref{hn}). Let
$$X=[2\kappa_2\kappa_{12}+\kappa_1(\kappa_{11}-\kappa_{22})]E_1+[2\kappa_1\kappa_{12}-\kappa_2(\kappa_{11}-\kappa_{22}]E_2.$$
Since
$\<X,Q\nabla\kappa\>=2D^2\kappa(\nabla\kappa,Q\nabla\kappa)=0,$ we
have
$$X=\frac{\<X,\nabla\kappa\>}{|\nabla\kappa|^2}\nabla\kappa.$$ A
simple computations shows that
$$\<X,\nabla\kappa\>=2D^2\kappa(\nabla\kappa,\nabla\kappa)-|\nabla\kappa|^2\Delta\kappa,$$
which completes the proof. \hfill$\Box$\\

It follows from Lemma \ref{l2} and Theorem \ref{2.1}  that

\begin{thm} \label{nt2.3}
If the Gaussian curvature function $\kappa$ is not constant on
$\Om,$ Then
$$\dim H^1_\kf(\Om,T)\leq1.$$
Moreover, if there is a point $o\in \Om$ such that
$$[\<Q\nabla\kappa,\nabla\Delta\kappa\>]^2+[D^2\kappa(\nabla\kappa,Q\nabla\kappa)]^2>0\quad\mbox{at}\quad
o,$$ then $H^1_\kf(\Om, T)=\{0\}.$
\end{thm}

\section{Infinitesimal Isometries}
\def\theequation{3.\arabic{equation}}
\hskip\parindent We shall give some characteristic conditions on a
function $w$ for which there exists a vector field $W$ such that
$V=W+wN$ is to be an infinitesimal isometry (Theorem \ref{th3.1}).

Let $M$ be a surface with the induced metric $g$ from $\R^3.$ Let
$N$ be the unit normal field of $M$. Let $\Pi$ be the second
fundamental form of $M$. Let $\Om\subset M$ be an open set. For
$w\in H^1_\is(\Om),$ there exists a unique vector field $W$ on
$\Om$, which is perpendicular to $H^1_\kf(\Om,T)$ in $H^1(\Om,T),$
such that $(W,w)$ is to be an infinitesimal isometry on $\Om$, that
satisfies\be DW(X,X)+w\Pi(X,X)=0\qfq X\in M_x,\,\,x\in
\Om,\label{3.1}\ee where $D$ is the Levi-Civita connection of the
induced metric $g$.

Let $o\in\Om$ be such that $\Om$ is star-shaped with respect to $o$.
Let the frame field $T$ and $E$ be given  in (\ref{n2.11}). Let
$$W=\var(t,\theta) T+\phi(t,\theta) E\qfq x=\F(t,\theta)\in\Om\cap\exp_o\Sigma(o),$$ where $\var=\<W,T\>$ and
$\phi=\<W, E\>.$  In the sequel all our computations are made on the
region $\Om\cap\exp_o\Sigma(o).$

Similar to (\ref{2.5}),  the relation (\ref{3.1})
is equivalent to \be \left\{\begin{array}{l}\var_t+w\Pi(T,T)=0,\\
f\phi_t-f_t\phi+\var_\theta+2fw\Pi(T,E)=0, \\
\phi_\theta+f_t\var+fw\Pi(E,E)=0,\\
\var(0)=\<W,\si(\theta)\>,\quad\phi(0)=\<W,\dot{\si}(\theta)\>.
\end{array}\right.\label{3.3}\ee

Let $\var$ solve  the first equation in (\ref{3.3}) with the initial
data $\var(0)=\<W,\si(\theta)\>.$ Then
\be\var=\<W_0,\si(\theta)\>-\int_0^tw\Pi_{11}ds.\label{dff}\ee As in
(\ref{2.8*}), a similar computation shows that $\phi$ solves the
second equation in (\ref{3.3}) if and only if it satisfies \be
\left\{\begin{array}{l}\phi_{tt}+\kappa\phi=P(w),\\
\phi(0)=\<W,\dot{\si}(\theta)\>,\end{array}\right.\label{3.3*}\ee
where  \be P(w)=-2w_1\Pi_{12}+w_2\Pi_{11}-w\Pi_{121}\qfq
x\in\Om,\label{3.4*}\ee
$$w_1=\<Dw,T\>,\quad w_2=\<DW,E\>,\quad
\Pi_{12}=\Pi(T,E),\quad\Pi_{121}=D\Pi(T,E,T),$$ etc. Furthermore,
differentiating the third equation in (\ref{3.3}) with respect to
the variable $t$ and using the first equation of (\ref{3.3}) yield
\be
0=\phi_{t\theta}+f_{tt}\var+f_t[w\Pi(E,E)-w\Pi(T,T)]+f[w\Pi(E,E)]_t\qfq
t>0.\label{c3.3}\ee Letting $t\rightarrow0$ in (\ref{c3.3}), we
obtain another initial data for the problem (\ref{3.3*})
\be\phi_t(0)=-w(o)\Pi(\si(\theta),\dot{\si}(\theta)).\label{cc3.3}\ee

Let $k$ be an integer. Let $T^k(M)$ be all tensor fields of rank $k$
on $M$.   Let
$${\bf R}_{XY}:\,\,T^k(M)\rightarrow T^k(M)$$ be the curvature
operator where $X,\,Y\in\X(M)$ are vector fields. For $K\in T^k(M),$
We have the following formulas, called the Ricci identities, \be
D^2K(\cdots, X,Y)=D^2K(\cdots, Y,X)+({\bf
R}_{XY}K)(\cdots).\label{3.5**}\ee The above formulas are very
useful when we have to exchange the order of the covariant
differentials of a tensor field.

 We seek some conditions on $w$ such that the
problem (\ref{3.1}) has a vector field solution $W$.

\begin{lem}\label{l3.1} Let $M$ be orientable. Let $(W,w)$ be an infinitesimal isometry of
$\Om$. Then\be \<D^2w,Q^*\Pi\>+w\kappa\tr\Pi=\<\nabla\kappa,W\>\qfq
x\in\Om,\label{3.2}\ee where $Q$ is defined by $(\ref{2.9*})$,
$\kappa$ is the Gaussian curvature function,  and $\nabla$, $\tr$
are the gradient, the trace of the induced metric of $M$,
respectively.
\end{lem}

{\bf Proof}\,\,\, Let $o\in \Om$ be any point. Then there is
$\varepsilon>0$ such that the geodesic ball centered at $o$ with the
radius $\varepsilon$ is contained in $\Om.$ Therefore, the systems
(\ref{3.3}) and (\ref{3.3*}) make sense for
$(t,\theta)\in[0,\varepsilon)\times[0,2\pi).$

From (\ref{3.3*}) and using the symmetry of $D\Pi$, we have
\beq&&\phi_{tt\theta}+\kappa_\theta\phi+\kappa\phi_\theta=-2(w_1\Pi_{12})_\theta+(w_2\Pi_{11})_\theta
-(w\Pi_{121})_\theta\nonumber\\
&&=-2f(w_{12}\Pi_{12}+w_1\Pi_{122})-2f_t[w_2\Pi_{12}+w_1(\Pi_{22}-\Pi_{11})]\nonumber\\
&&\quad+f(w_{22}\Pi_{11}+w_2\Pi_{112})+f_t(2w_2\Pi_{12}-w_1\Pi_{11})\nonumber\\
&&\quad-f(w_2\Pi_{121}+w\Pi_{1212})-f_tw(2\Pi_{221}-\Pi_{111})\nonumber\\
&&=f(w_{22}\Pi_{11}-2w_{12}\Pi_{12}-2w_1\Pi_{122}-w\Pi_{1212})\nonumber\\
&&\quad+f_t[w_1(\Pi_{11}-2\Pi_{22})+w(\Pi_{111}-2\Pi_{221})],\label{3.6*}\eeq
which yields \beq
\phi^{(3)}_\theta(0)&&=w_{22}\Pi_{11}-2w_{12}\Pi_{12}-2w_1\Pi_{122}-w\Pi_{1212}\nonumber\\
&&\quad+[w_1(\Pi_{11}-2\Pi_{22})+w(\Pi_{111}-2\Pi_{221})]'(0)\nonumber\\
&&\quad-\kappa_\theta'\phi(0)-\kappa'\phi_\theta(0)-\kappa\phi'_\theta(0)\nonumber\\
&&=w_{11}(\Pi_{11}-2\Pi_{22})+w_{22}\Pi_{11}-2w_{12}\Pi_{12}+2w_1(\Pi_{111}-3\Pi_{221})\nonumber\\
&&\quad+w[\Pi_{1111}-3\Pi_{2211}+2\kappa(\Pi_{22}-\Pi_{11})]\nonumber\\
&&\quad-\<\nabla\kappa,\dot{\si}(\theta)\>\<W,\dot{\si}(\theta)\>+\<\nabla\kappa,\si(\theta)\>\<W,\si(\theta)\>,
\label{3.5}\eeq  where the following formulas have been used
$$\Pi_{1212}=\Pi_{2211}+{\bf R}_{TE}D^2\Pi(T,E)=\Pi_{2211}+\kappa(\Pi_{11}-\Pi_{22})\,(\mbox{by}\,(\ref{3.5**})),$$
$$\phi'_\theta(0)=w\kappa(\Pi_{11}-\Pi_{22}).$$

On the other hand, using the equation (\ref{2.2}) and the first
equation in (\ref{3.3}), we obtain\beq
(f_t\var)^{(3)}(0)&&=[f^{(4)}\var+3f^{(3)}\var'+3f''\var''+f'\var^{(3)}](0)\nonumber\\
&&=-2\kappa'\var(0)-3\kappa(0)\var'(0)+\var^{(3)}(0)\nonumber\\
&&=-2\<\nabla\kappa,\si(\theta)\>\<W,\si(\theta)\>
-w_{11}\Pi_{11}-2w_1\Pi_{111}\nonumber\\
&&\quad+w(3\kappa\Pi_{11}-\Pi_{1111})\quad\mbox{at}\quad
o.\label{3.6}\eeq  Moreover, we have \beq
(fw\Pi_{22})^{(3)}(0)&&=f^{(3)}(0)w(o)\Pi_{22}(o)+3(w\Pi_{22})''(0)\nonumber\\
&&=3w_{11}\Pi_{22}+6w_1\Pi_{221}+w(3\Pi_{2211}-\kappa
\Pi_{22})\quad\mbox{at}\quad o.\label{3.7}\eeq

Finally, using the third equation in (\ref{3.3}), we obtain from
(\ref{3.5})-(\ref{3.7})\beq
0&&=(\phi_\theta+f_t\var+fw\Pi_{22})^{(3)}(0)\nonumber\\
&&=w_{11}\Pi_{22}-2w_{12}\Pi_{12}+w_{22}\Pi_{11}+w\kappa(\Pi_{11}+\Pi_{22})\nonumber\\
&&\quad-\<\nabla\kappa,\dot{\si}(\theta)\>\<W,\dot{\si}(\theta)\>-\<\nabla\kappa,\si(\theta)\>\<W,\si(\theta)\>\nonumber\\
&&\quad=\<D^2w,Q^*\Pi\>+w\kappa\tr\Pi-\<\nabla\kappa,W\>\quad\mbox{at}\quad
o.\nonumber\eeq   \hfill$\Box$

Let $s\geq0$ be given. Let $\Phi_0(t)$ and $\Phi(t,s)$ solve the
problem \be \left\{\begin{array}{l}
\Phi_{0tt}(t)+\kappa(t)\Phi_0(t)=0\qfq t\geq 0,\\
\Phi_0(0)=1,\quad \Phi_{0t}(0)=0,\end{array}\right.\label{3.8}\ee
and \be \left\{\begin{array}{l}
\Phi_{tt}(t,s)+\kappa(t)\Phi(t,s)=0\qfq t\geq s,\\ \Phi(s,s)=0,\quad
\Phi_t(s,s)=1,\end{array}\right.\label{3.8*}\ee respectively. Note
that
$$\Phi(t,0)=f.$$
Let $w$ be a function on $\Om$ and $W_o\in M_o$. Let
\be\phi=\Phi_0(t)\<W_o,\dot{\si}(\theta)\>-w(o)\Pi(\dot{\si}(\theta),\si(\theta))f
+\int_0^t\Phi(t,s)P(w)(s)ds,\label{3.11}\ee where $P(w)$ is given by
(\ref{3.4*}). Then $\phi$ solves the problem (\ref{3.3*}) and
(\ref{cc3.3}).

We have

\begin{thm} \label{th3.1}Let $M$ be orientable and let $\Om$ be star-shaped with respect to a point $o\in\Om$.
Then $w\in H^1_\is(\Om)$ if and only if it is in the form of\be
w=u(x)+\<W_o,N\>(x)\qfq x\in\Om\cap\exp_o\Sigma(o),\label{3.12f}\ee
where $W_o\in M_o$ is a constant vector and $u$ is a solution to the
problem \be \A_o
u+u(o)\Pi(\si(\theta),\dot{\si}(\theta))\kappa_2f=0\qfq
x\in\Om\cap\exp_o\Sigma(o),\label{3.12f2}\ee where\be \A_o
u=\<D^2u,Q^*\Pi\>+u\kappa\tr\Pi+\kappa_1\int_0^tu\Pi_{11}ds
-\kappa_2\int_0^t\Phi(t,s)P(u)(s)ds.\label{3.12}\ee
\end{thm}

\begin{rem}Since $\Om$ is star-shaped with respect to $o$, a
simply computation shows that \be |u(o)|\leq
C(\|u\|_{H^1(\Om)}+\|u\|_{L^2(\Ga)})\qfq u\in H^1(\Om).\ee Then
the second term in the left hand side of $(\ref{3.12f2})$ makes
sense for $u\in H^1(\Om).$
\end{rem}

\begin{rem}
As a constant vector $M_o$ on $\Om$, or a translation displacement
of $\Om$, $(\tilde{W}_o,\<W_o N\>)$ is a trivial  infinitesimal
isometry where $W_o=\tilde{W}_o+\<W_o,N\>N.$ Then a solution $u$ to
the problem $(\ref{3.12f2})$ is itself in $H^1_\is(\Om).$
\end{rem}

\begin{rem} The formula $(\ref{3.12f})$ depends on the choice of the point $o\in\Om.$
If the point $o$ can be chosen to be
an umbilical point of $M$ $(\cite{dC})$, then $\kappa(o)\geq0$ and
$$\Pi(o)=\sqrt{\kappa(o)}g,$$ which yields
$$\Pi(o)(\si(\theta),\dot{\si}(\theta))=0\qfq \theta\in(0,2\pi].$$
In this case the equation $(\ref{3.12f2})$ becomes \be \A_ou=0\qfq
x\in\Om.\label{3.12f3}\ee

Another case for which $(\ref{3.12f3})$ holds true is that $o\in\Om$
can be chosen such that \be\kappa_2=0\qfq x\in\Om.\label{3.12f4}\ee
\end{rem}

{\bf Proof of Theorem 4.1}\,\,\,\,\,{\bf Necessity}\,\,\,Let $w\in
H^1_\is(\Om).$ Let a vector field $W\bot H^1_\kf(\Om,T)$ be such
that $(W,w)$ is an infinitesimal isometry. Let
$$W(o)=\hat W(o)+\<W(o),N\>N\qfq x\in\Om.$$
Let
$$U=W-\hat W(o),\quad u=w-\<W(o),N\>\qfq x\in\Om.$$ Then
$(U,u)$ is  an infinitesimal isometry field with $U(o)=0.$ Using the
formulas (\ref{dff}) and (\ref{3.11}) in the formula (\ref{3.2})
where $(W,w)$ is replaced by $(U,u)$, we have the formula
(\ref{3.12}) for $u$.

{\bf Sufficiency}\,\,\,Let $u$ solve the problem (\ref{3.12f2}). It
will suffice to prove that there is a vector field $U\in\X(\Om)$
such that $(U,u)$ is an infinitesimal isometry. We define
$$U=\var T+\phi E\qfq x\in\Om\cap\exp_o\Sigma(o),$$ where
\be\var=-\int_0^tu\Pi_{11}ds,\quad\phi=-u(o)\Pi(\dot{\si}(\theta),\si(\theta))f+\int_0^t\Phi(t,s)P(u)ds.\label{dy}\ee
Then the equation (\ref{3.12f2}) means that \be
\<D^2u,Q^*\Pi\>+u\kappa\tr\Pi=\<\nabla\kappa,U\>\qfq
x\in\Om\cap\exp_o\Sigma(o).\label{3.13}\ee

Clearly, $\var$ and $\phi$ satisfy the first equation and the second
equation in (\ref{3.3}).  To complete the proof, it remains to show
that $\var$ and $\phi$, given by (\ref{dy}), solve the third
equation  in (\ref{3.3}). For this end, we let
$$\eta=\phi_\theta+f_t\var+fu\Pi_{22}\qfq x\in\Om\cap\exp_o\Sigma(o).$$
 Using (\ref{2.2}), (\ref{3.6*}), and
(\ref{3.13}),  we compute \beq
\eta''&&=\phi''_\theta+f^{(3)}\var+2f''\var'+f'\var''+f(u\Pi_{22})''+2f'(u\Pi_{22})'+f''u\Pi_{22}\nonumber\\
&&=\phi''_\theta-(f\kappa'+f'\kappa)\var-2f\kappa\var'+f'\var''+f(u\Pi_{22})''+2f'(u\Pi_{22})'-f\kappa
u\Pi_{22}\nonumber\\
&&=\phi''_\theta+f[(u\Pi_{22})''-\kappa
u\Pi_{22}-2\kappa\var'-\kappa'\var]+f'[2(u\Pi_{22})'-\kappa\var+\var'']\nonumber\\
&&=-(\kappa_\theta\phi+\kappa\phi_\theta+f\kappa'\var+f'\kappa\var)\nonumber\\
&&\quad+f[(u\Pi_{22})''-\kappa
u\Pi_{22}-2\kappa\var']+f'[2(u\Pi_{22})'+\var'']\nonumber\\
&&\quad+f[u_{22}\Pi_{11}-2u_{12}\Pi_{12}-2u_1\Pi_{122}-u\Pi_{1122}]\nonumber\\
&&\quad+f'[u_1(\Pi_{11}-2\Pi_{22})+u(\Pi_{111}-2\Pi_{122})]\nonumber\\
&&=-[f\<\nabla\kappa,U\>+\kappa(\phi_\theta+f'\var+fu\Pi_{22})]+f(\<D^2u,Q^*\Pi\>+u\kappa\tr\Pi)\nonumber\\
&&\quad+f'(u_1\Pi_{11}+u\Pi_{111}+\var'')\nonumber\\
&&=f(\<D^2u,Q^*\Pi\>+u\kappa\tr\Pi-\<\nabla\kappa,U\>)-\kappa\eta+f'(u\Pi_{11}+\var')'\nonumber\\
&&=-\kappa\eta,\label{3.15}\eeq where the following formula has been
used
$${\bf R}_{E_1E_2}\Pi(E_1,E_2)=\kappa(\Pi_{11}-\Pi_{22})(\mbox{by}\,(\ref{3.5**})).$$
Moreover, we have the initial data
$$\eta(0)=\phi_\theta(0)+\var(0)=0,\quad\eta'(0)=\phi'_\theta(0)+\var'(0)+u(o)\Pi_{22}(o)=0,$$
which imply by the equation (\ref{3.15}) that $(U,u)$ is an
infinitesimal isometry.  \hfill$\Box$ \\

If  surface $M$ is given as a graph, an infinitesimal isometry
function $w\in H^1_\is(\Om)$ can be written as an explicit formula
in the Cartesian orthogonal coordinate system. Let \be
M=\{\,(x,h(x))\,|\,x=(x_1,x_2)\in\R^2\,\},\label{4.17}\ee where $h$
is a smooth function on $\R^2.$  Let \be V(p)=(u_1,u_2,u)\qfq p\in
M.\label{3.19*}\ee

We then have

\begin{thm}$(\cite{Al}, \cite{Pog1}, \cite{Pog}, \cite{Sp})$ \label{t4.2}
Let $\tilde{\Om}\subset\R^2$ be a star-shaped with respect to a
point $\tilde{o}\in\tilde{\Om}$. Then there are functions $u_1,$
$u_2$ such that $V$ is an infinitesimal isometry  on
$$\Om=\{\,(x,h(x))\,|\,x\in\tilde{\Om}\,\}$$ if and only if $u$ solves the problem \be
\tilde{\div}A(x)\tilde{\nabla} u=0\qfq
x\in\tilde{\Om},\label{3.18}\ee where $\tilde{\div}$ and
$\tilde{\nabla}$ are the divergence and gradient of $\R^2$ in the
Euclidean metric, respectively, and \be
A(x)=\left(\begin{array}{cc}h_{x_2x_2}&-h_{x_1x_2}\\
-h_{x_1x_2}&h_{x_1x_1}\end{array}\right)\qfq x\in\tilde{\Om}.\ee
\end{thm}

Based on Theorem \ref{t4.2}, we shall work out a formula for $w$ as
follows.

It follows from (\ref{4.17}) that
$$M_p=\{\,(\a,\a(h))\,|\,\a=(\a_1,\a_2)\in\R^2\,\}\qfq
p=(x,h(x))\in M.$$ Moreover,
$$N=\eta(\t{\nabla} h,-1),\quad \eta=\frac{1}{\sqrt{1+|\t{\nabla} h|^2}},$$
\be\Pi((\a,\a(h)),(\b,\b(h)))=\eta\D^2h(\a,\b)\qfq
(\a,\a(h)),\,\,(\b,\b(h))\in M_p,\label{4.8}\ee
$$\kappa(p)=\eta^4(h_{x_1x_1}h_{x_2x_2}-h^2_{x_1x_2})\qfq p\in M,$$ where  $\D^2h$ is the Hessian of $h$ in the
Euclidean metric of $\R^2.$

Let $o=(0,0,0).$ Consider the polar coordinates $(t,\theta)$ in the
induced metric $g$ of $M$ initiating from $o.$  Consider a
1-parameter family of geodesics given by
$$ {\cal F}(t,\theta)=\exp_ot\si(\theta),$$ where
 \be\si(\theta)=\frac{\cos\theta}{\sqrt{1+h_{x_1}^2}}(1,0,h_{x_1})+\frac{\sin\theta}{\sqrt{(1+h_{x_1}^2)(1+|\D h|^2)}}
 (-h_{x_1}h_{x_2},1+h_{x_1}^2,h_{x_2}).\label{the}\ee
 Let $$ {\cal
F}(t,\theta)=(r(t),h(r(t))),$$ where $r(t)=(r_1(t),r_2(t))$ is a
curve in $\R^2.$ Then $$\Big(\ddot{r},\,\D^2h(\dot{r},\dot{r})+\<\D
h,\ddot{r}\>\Big)=\D_{\dot{\F}}\dot{\F}=D_{\dot{\F}}\dot{\F}-\Pi(\dot{\F},\dot{\F})N=
-\Pi(\dot{r},\dot{r})N,$$ which yields \be
\ddot{r}(t)+\eta^2\D^2h(\dot{r}(t),\dot{r}(t))\D h=0\qfq
t>0,\label{the1}\ee with the initial data \be
r(0)=0,\quad\dot{r}(0)=\si(\theta). \label{the2}\ee

Set \be X_i(t,s)=\D_{\dot{r}(s)}[(r_i(s)-r_i(t))Q\t{\nabla}h(s)]\qfq
x=r(t)\in R^2,\,\,i=1,\,\,2,\label{3.23*}\ee where
$Q\t{\nabla}h=(h_{x_2},-h_{x_1}).$

\begin{thm} Let $M$ be given by $(\ref{4.17})$ and let $w\in H^1_\is(\Om)$ be given.
Then there is a solution $u$ to the problem $(\ref{3.18})$ such
that\beq w/\eta&&=\<Z,\D
h\>+h_{x_1}(x)\int_0^t[X_2(t,s)(u)-u_{x_1}\dot{r}(h)]ds\nonumber\\
&&\quad
-h_{x_2}(x)\int_0^t[X_1(t,s)(u)+u_{x_2}\dot{r}(h)]ds-u(x)\qfq
x=r(t)\in R^2,\label{3.23}\eeq where $Z\in\R^2$ is a constant
vector.
\end{thm}

{\bf Proof}\,\,\,Since $w\in H^1_\is(\Om),$ there is a unique
displacement field $V=(u_1,u_2,u)$ in $\R^3$ such that
$$w=\<V,N\>,\quad\<\D_XV,X\>=0\qfq X\in M_x,\,\,x\in\Om.$$
By Theorem \ref{t4.2}, the third component $u$ of $V$ solves the
problem (\ref{3.18}).   We shall obtain $u_1$ and $u_2$ in terms of
$u$.

By Lemma 3 of Chapter 12 in \cite{Sp}, there is a vector field
$Y=(\psi_1,\psi_2,\psi)$ such that \be \D_XV=X\times Y\qfq X\in
M_x,\,\,x\in\Om,\label{3.24}\ee where $\times$ denotes the exterior
product.  Letting $X=(1,0,h_{x_1})$ and $X=(0,1,h_{x_2})$ in
(\ref{3.24}), respectively, we obtain \be
Y=(-u_{x_2},u_{x_1},\psi),\ee \be u_{1x_1}=-h_{x_1}u_{x_1},\quad
u_{1x_2}=\psi-h_{x_2}u_{x_1},\quad
u_{2x_1}=-\psi-h_{x_1}u_{x_2},\quad
u_{2x_2}=-h_{x_2}u_{x_2}.\label{3.26}\ee It follows from
(\ref{3.26}) that \be
\left\{\begin{array}{l}u_1=z_1+\int_0^t[\psi\dot{r}_2-u_{x_1}\dot{r}(h)]ds,\\
u_2=z_2-\int_0^t[\psi\dot{r}_1+u_{x_2}\dot{r}(h)]ds,\end{array}\right.\label{3.27}\ee
where $z_1$ and $z_2$ are constants.

Next, we compute $\psi.$ Since the curvature operator of $\R^3$ in
the Euclidean metric is  zero, \be
-\D_X\D_ZV+\D_Z\D_XV+\D_{[X,Z]}V=0\label{3.25}\ee for vector fields
$X,$ $Z$ on $\Om.$

Using the formula (\ref{3.24}) in the identity (\ref{3.25}), we
obtain $$X\times\D_ZY=Z\times\D_XY.$$ In particular, taking
$X=(1,0,h_{x_1})$ and $Z=(0,1,h_{x_2})$, respectively,  yield
$$\left\{\begin{array}{l}\psi_{x_1}=-h_{x_1}u_{x_1x_2}+h_{x_2}u_{x_1x_1},\\
\psi_{x_2}=-h_{x_1}u_{x_2x_2}+h_{x_2}u_{x_1x_2},
\end{array}\right.$$  which give
\be
\psi=\<\t{\nabla}u,Q\t{\nabla}h\>-\int_0^t\<\t{\nabla}u,\D_{\dot{r}}Q\t{\nabla}h\>ds.\label{3.29}\ee
Inserting the formula (\ref{3.29}) into the formula (\ref{3.27}), we
have\be\left\{\begin{array}{l}u_1=z_1+\int_0^t[X_2(t,s)(u)-u_{x_1}\dot{r}(h)]ds,\\
u_2=z_2-\int_0^t[X_1(t,s)(u)+u_{x_2}\dot{r}(h)]ds.\end{array}\right.\label{3.30}\ee
Then the formula (\ref{3.23}) follows from (\ref{3.30}).    \hfill$\Box$\\

\section{Elliptic Surfaces}
\def\theequation{4.\arabic{equation}}
\hskip\parindent  Let $M$ be a surface in $\R^3$.  $M$ is said to be
elliptic if the fundamental form $\Pi$ is positive for all $x\in M.$
Assume that $M$ be elliptic throughout this section. Then the
problem (\ref{3.12f2}) will become an elliptic one (Theorem
\ref{t5.1}).

We introduce another metric on $M$ by
$$\hat{g}=\Pi\qfq x\in M.$$

\begin{pro}\label{p4.1}  Let $M$ be elliptic. Then for $w\in C^2(M)$,
\be \kappa\Delta_\Pi w+\frac{1}{2\kappa}Q^*\Pi(\nabla\kappa,\nabla
w)=\<D^2w,Q^*\Pi\>\qfq x\in M,\label{4.1}\ee where $\Delta_\Pi$ is
the Laplacian of the metric $\hat{g}=\Pi$ and
$Q:\X(\Om)\rightarrow\X(\Om)$ is the operator, given by
$(\ref{2.9*}).$
\end{pro}

{\bf Proof}\,\,\, Let $o\in M$ be fixed. Consider the polar
coordinates in the induced metric $g$
$$\pl t=T,\quad\pl\theta=fE.$$ Note that the above
$(\pl t,\pl\theta)$ is no longer the polar coordinates in the metric
$\hat{g}=\Pi.$

In the coordinate system $(\pl t,\pl\theta)$, we have
$$\hat{g}=\hat{g}_{11}dt^2+\hat{g}_{12}(dt
d\theta+d\theta dt)+\hat{g}_{22}d\theta^2,$$
$$\hat{G}=\Big(\hat{g}_{ij}\Big)=\left(\begin{array}{cc}
\Pi_{11}&f\Pi_{12}\\
f\Pi_{12}&f^2\Pi_{22}\end{array}\right),\quad \det \hat{G}=\kappa
f^2,\quad \hat{G}^{-1}=\frac{1}{\kappa f^2}\left(\begin{array}{cc}
f^2\Pi_{22}&-f\Pi_{12}\\
-f\Pi_{12}&\Pi_{11}\end{array}\right).$$ Moreover,
$$w_{t\theta}=fw_{12}+f'w_2,\quad w_{\theta\theta}=f^2w_{22}-ff'w_1+f_\theta w_2.$$
Using those formulas, we obtain \beq\kappa\Delta_\Pi
w&&=\frac{\kappa}{\sqrt{\kappa}f}[(\sqrt{\kappa}f\frac{\Pi_{22}}{\kappa}w_t)_t
-(\sqrt{\kappa}f\frac{\Pi_{12}}{\kappa
f}w_t)_\theta-(\sqrt{\kappa}f\frac{\Pi_{12}}{\kappa
f}w_\theta)_t+(\sqrt{\kappa}f\frac{\Pi_{11}}{\kappa
f^2}w_\theta)_\theta]\nonumber\\
&&=\<D^2w,Q^*\Pi\>+\{\frac{\sqrt{\kappa}}{f}[(\frac{f\Pi_{22}}{\sqrt{\kappa}})_t-(\frac{\Pi_{12}}{\sqrt{\kappa}})_\theta]
-\Pi_{11}\}w_1-2\frac{f'}{f}\Pi_{12}w_2\nonumber\\
&&\quad+\{\sqrt{\kappa}[(\frac{\Pi_{11}}{\sqrt{\kappa}f})_\theta-(\frac{\Pi_{12}}{\sqrt{\kappa}})_t]
+\frac{f_\theta}{f^2}\Pi_{11}\}w_2\nonumber\\
&&=\<D^2w,Q^*\Pi\>+\frac{1}{2\kappa}(\kappa_2\Pi_{12}-\kappa_1\Pi_{22})w_1-2\frac{f'}{f}\Pi_{12}w_2\nonumber\\
&&\quad+[2\frac{f'}{f}\Pi_{12}+\frac{1}{2\kappa}
(\kappa_1\Pi_{12}-\kappa_2\Pi_{11})]w_2\nonumber\\
&&=\<D^2w,Q^*\Pi\>-\frac{1}{2\kappa}\Pi(Q\nabla\kappa,Q\nabla
w).\nonumber\eeq      \hfill$\Box$

Let $\Om\subset M$ with boundary $\Ga.$  Instead of the usual inner
product of $L^2(\Om)$, we use the following inner product on
$L^2(\Om)$
$$(w,v)_{L^2_\Pi(\Om)}=\int_\Om wvdg_\Pi\qfq w,\,\,v\in L^2(\Om).$$
We denote by $L^2_\Pi(\Om)$ the above space.

It is well known that the negative Laplacian operator $-\Delta_\Pi$
on $\Om$ with the Dirichlet boundary condition is a positive
selfadjoint operator on $L^2_\Pi(\Om)$ and
$$D(\Delta_\Pi)=H^2(\Om)\cap H^1_0(\Om).$$  Moreover, we  extend the domain $H^2(\Om)\cap H_0^1(\Om)$
of $\Delta_\Pi$ to $H^1_0(\Om)$ such that\be \Delta_\Pi:\quad
H^1_0(\Om) \rightarrow H^{-1}(\Om)\label{tg}\ee is an isomorphism.
We let \be \B
w=\B_ow+w(o)\frac{\kappa_2}{\kappa}\Pi(\si(\theta),\dot{\si}(\theta))f\qfq
x\in\Om\cap\exp_o\Sigma(o), \label{n4.3}\ee where \beq \B_o
w&&=\frac{1}{2\kappa^2}Q^*\Pi(\nabla\kappa,\nabla
w)+w\tr\Pi+\frac{\kappa_1}{\kappa}\int_0^tw\Pi_{11}ds\nonumber\\
&&\quad -\frac{\kappa_2}{\kappa}\int_0^t\Phi(t,s)P(w)(s)ds \qfq
x\in\Om\cap\exp_o\Sigma(o).\label{nn4.3}\eeq

\begin{rem} Let $o\in\Om.$ Since $\Om\cap\exp_oC(o)$ is zero
mensurable and
$$\Om=[\Om\cap\exp_o\Sigma(o)]\cup[\Om\cap\exp_oC(o)],$$
$\B w$ is defined by $(\ref{n4.3})$ on $\Om$ almost everywhere.
\end{rem}

Consider the operator $\A_o$, defined by (\ref{3.12}). It follows
from (\ref{3.12}) and (\ref{4.1}) that \be \A_o
w+w(o)\kappa_2\Pi(\si(\theta),\dot{\si}(\theta))=\kappa(\Delta_\Pi
w+\B w).\label{fl}\ee

Since for $W_o\in M_o$, $(\tilde{W}_o,\<W_o,N\>)$ is a trivial,
smooth infinitesimal isometry where $W_0=\tilde{W}_o +\<W_o,N\>N,$
it follows from Theorem \ref{th3.1} and (\ref{fl}) that

\begin{thm}\label{t5.1} Let $\Om\subset M$ be elliptic and star-shaped  with
respect to $o\in\Om.$ Then \be H^1_\is(\Om)=\{\,\,w\,\,|\,\,w\in
H^1(\Om),\,\,\Delta_\Pi w+\B w=0\,\,\}.\ee
\end{thm}

Next, we consider the structure of solutions to the equation
$\Delta_\Pi w+\B w=0$ in $H^1(\Om).$

Since $\Delta_\Pi^{-1}w\in H_0^1(\Om)$ for $w\in H^{-1}(\Om),$ we
have the following estimates
$$\|\B_o\Delta^{-1}_\Pi w\|_{H^{-1}(\Om)}\leq C\|\Delta^{-1}_\Pi
w\|_{L^2_\Pi(\Om)}\leq C\|\Delta^{-1}_\Pi w\|_{H^1_0(\Om)}\leq
C\|w\|_{H^{-1}(\Om)}\qfq w\in L^2(\Om),$$ which yield

\begin{lem}\label{l5.1} The operator $\B_o\Delta_\Pi^{-1}:\,\,\,H^{-1}(\Om)\rightarrow H^{-1}(\Om)$ is a compact
operator.
\end{lem}

Consider the operator $\Delta_\Pi+\B$ with the domain
$D(\Delta_\Pi+\B)= H^1_0(\Om)$. Denote by $\B^*$ the adjoint
operator of $\B$ with respect to the inner product of
$L^2_\Pi(\Om).$ Then \be\B^*=\B_o^*+
\Big(\frac{\kappa_2}{\kappa}[\Pi(\si(\theta),\dot{\si}(\theta))f],\quad\cdot\Big)_{L^2_\Pi(\Om)}\delta(o),\label{DR}\ee
where $\delta(o)\in H^{-1}(\Om)$ is the Dirac function at $o$ and
$$D(\Delta_\Pi+\B^*)= H^1_0(\Om).$$ Let \be
{{\V_0}}(\Om)=\{\,\var\,|\,\var\in
H^1_0(\Om),\,\,\Delta_\Pi\var+\B\var=0\,\},\label{58}\ee
 \be
{\V_0}_*(\Om)=\{\,\var\,|\,\var\in
H^1_0(\Om),\,\,\Delta_\Pi\var+\B^*\var=0\,\},\ee \be
{\V_0}_*(\Ga)=\{\,\var_\nu\,|\,\var\in{\V_0}_*(\Om)\,\}.\ee It
follows from Lemma \ref{l5.1} and the formula (\ref{DR}) that
$$\Delta_\Pi^{-1}\B^*=(\B_o\Delta_\Pi^{-1})^*+
\Big(\frac{\kappa_2}{\kappa}[\Pi(\si(\theta),\dot{\si}(\theta))f],\quad\cdot\Big)_{L^2_\Pi(\Om)}\Delta_\Pi^{-1}\delta(o):$$
$L^2_\Pi(\Om)\rightarrow L^2_\Pi(\Om)$ is a compact operator. Then
$\V_0(\Om)$ and $\V_0(\Ga)$ are subspaces of finite dimension.

We discompose $H^1_\is(\Om)$  as a direct sum in $H^1(\Om)$ as \be
H^1_\is(\Om)=\V_0(\Om)\oplus \V_0^\bot(\Om).\label{nn4.11}\ee

We have
\begin{thm}\label{t5.2} Let ${\V_0}_*(\Om)=\{\,0\,\}.$
Then $w\in\V_0^\bot(\Om)$ if and only if there is a unique $\psi\in
H^{1/2}(\Ga)$ such that \be
w=w_0-\Delta^{-1}_\Pi(I+\B\Delta^{-1}_\Pi)^{-1}\B
w_0,\label{3.19*}\ee where $w_0\in H^1(\Om)$ is the unique solution
to the problem \be \left\{\begin{array}{l}\Delta_\Pi w_0=0\qfq
x\in\Om,\\
 w_0=\psi\qfq x\in\Gamma.\end{array}\right.\label{3.20*}\ee

If ${\V_0}_*(\Om)\not=\{\,0\,\},$ then $w\in\V_0^\bot(\Om)$ if and
only if there is a unique $\psi\in H^{1/2}(\Ga)$
satisfying\be(\psi,\var_\nu)_{L^2_\Pi(\Ga)}=0\quad\mbox{for
all}\quad \var\in {\V_0}_*(\Om),\label{3.11**}\ee  such that
$(\ref{3.20*})$ and $(\ref{3.19*})$ hold.
\end{thm}

{\bf Proof}\,\,\, By Theorem \ref{th3.1},   what we are looking for
is a solution $w\in H^1_\is(\Om)$ to the problem \be
\left\{\begin{array}{l}\Delta_\Pi
w+\B w=0\qfq x\in\Om,\\
w=\psi\qfq x\in\Ga.\end{array}\right.\label{5.9}\ee Let $w_0\in
H^1(\Om)$ be the solution to the problem (\ref{3.20*}) and let
$v=w-w_0.$ Then the problem (\ref{5.9}) is equivalent to solve
 \be\Delta_\Pi v+\B v=-\B w_0\quad\mbox{for some}\quad v\in
 H^1_0(\Om).\label{5.10}\ee Let $u=\Delta_\Pi v.$ Then
 the problem (\ref{5.10}) is the same to the problem \be
 u+\B\Delta^{-1}_\Pi u=-\B w_0.\label{5.11}\ee By the Fredholm theorem (\cite{Kr}), the problem
 (\ref{5.11}) is solvable if and only if \be
 (\B w_0,\,\,\var)_{L^2_\Pi(\Om)}=0\label{5.12}\ee for all
  $\var\in\V$ where
 \be \V=\{\,\var\in L^2_\Pi(\Om)\,|\,\var+(\B\Delta^{-1}_\Pi)^*\var=0\,\}.\label{5.13}\ee
It is easy to check that $$\V={\V_0}_*(\Om)=\{\,\var\in
H^1_0(\Om)\,|\,\Delta_\Pi\var+\B^*\var=0\,\}.$$

Clearly, if ${\V_0}_*(\Om)=\{\,0\,\},$ the claim is true. We assume
that ${\V_0}_*(\Om)\not=\{\,0\,\}.$   Let $\psi\in H^{1/2}(\Ga)$ be
such that (\ref{3.11**}) is true. Then it follows from
(\ref{3.11**}) that \beq&&(\B
w_0,\,\,\var)_{L^2_\Pi(\Om)}=(w_0,\B^*\var)_{L^2_\Pi(\Om)}
=-(w_0,\Delta_\Pi\var)_{L^2_\Pi(\Om)}=-(\psi,\var_\nu)_{L^2_\Pi(\Ga)}=0,\nonumber\eeq
for all $\var\in{\V_0}_*(\Om).$ Then the problem (\ref{5.9}) has a
solution $w$, given by (\ref{3.19*}). Moreover, a simply computation
shows that if $w$ is a solution to the problem (\ref{5.9}), then the
conditions (\ref{3.11**}) hold. Then the proof is complete.
\hfill$\Box$\\

If $\Om$ is of constant curvature, it follows from the formulas
(\ref{n4.3}) and (\ref{nn4.3}) that
$$\B w=\frac{1}{2\kappa^2}Q^*\Pi(\nabla\kappa,\nabla
w)+w\tr\Pi.$$ We then have $\B w\in C^\infty(\Om)$ whenever $w\in
C^\infty(\Om).$ Let $\Om$ be star-shaped respect to $o$ and be not
of constant curvature. Consider a solution $w_0$ to the problem
(\ref{3.20*}). If $\psi\in C^\infty(\Ga)$, then $w_0\in
C^\infty(\Om)$. Furthermore, by the formula (\ref{n4.3}), $\B w_0$
is $C^\infty$ on $\Om\cap\exp_o\Sigma(o).$ Therefore from the
formula (\ref{3.19*}), $w$ is also $C^\infty$ on
$\Om\cap\exp_o\Sigma(o).$ Since $C^\infty(\Ga)$ is dense in
$H^{1/2}(\Ga)$, the ellipticity of the operator $\Delta_\Pi$ implies
the following density result.

\begin{thm}\label{tn4.3} Let $\Om\subset M$ be elliptic which is star-shaped
with respect to $o\in\Om$. Moreover, suppose that one of the
following assumptions holds true: $\Om$ is of constant curvature, or
\be\Om\subset\exp_o\Sigma(o).\label{n4.20}\ee Then the strong
$H^1(\Om)$ closure of
$$H^1_{\ib}(\Om)\cap C^\infty(\Om)$$ agrees with $H^1_{\ib}(\Om).$
\end{thm}

\begin{rem} Let $\Om$ be bounded and be not of constant curvature.
If $\Om\cap\exp_oC(o)\not=\emptyset,$ we only have
$$\B u\in L^2(\Om)$$ when $u\in C^\infty(\Om)$ where the operator $\B$ is defined by
$(\ref{n4.3}).$  A condition like $(\ref{n4.20})$ is necessary for
the above density result. An interesting example is given by
$\cite{CoVo}$ $($see $\cite{LeMoPa2}$ or $\cite{Sp}$$)$, where $\Om$
is a closed smooth surface of non-negative curvature for which
$C^\infty$ infinitesimal isometries consist only of trivial fields,
whereas there exist  non-trivial $C^2$ infinitesimal isometries.
Therefore $H^1_\is(\Om)\cap C^\infty(\Om)$ is not dense in
$H^1_\is(\Om)$ for this surface.
\end{rem}

By Theorem \ref{t5.2}, if $\V_0(\Om)=\{\,0\,\},$ an infinitesimal
isometry function $w\in H^1_\is(\Om)$ is completely given by its
boundary trace $w\in H^{1/2}(\Ga).$ However, in general
$\V_0(\Om)\not=\{\,0\,\}$ even for a spherical cap,  see Theorem
\ref{t4.6} later. Next, we consider several cases for which the
relations $\V_0(\Om)=\{\,0\,\}$ hold. This problem closely relates
to the first eigenvalue of $-\Delta_\Pi$ with the Dirichlet boundary
condition. Let $\lam_1$ be the first positive eigenvalue of
$-\Delta_\Pi$ on $L^2_\Pi(\Om).$

Let $\kappa_\Pi$ be the  curvature function of $M$ in the metric
$g_\Pi=\Pi.$  Let $\rho_\Pi=\rho_\Pi(x,o)$ be the distance
function from $x\in M$ to $o\in M$ in the metric $g_\Pi=\Pi.$ For
$a>0,$ let
$$\mu(a)=\sup_{\rho_\Pi\leq a}\kappa_\Pi.$$
Then $\mu(a)$ is an increasing function in $a\in[0,\infty).$ Let
$a_0>0$ be given by
$$\mu(a_0)a^2_0=\frac{\pi}{2}.$$

\begin{lem} Assume that there is $0<a<a_0$ such that $$\Om\subset\{\,x\,|\,x\in
M,\,\,\rho_\Pi(x)<a\,\}.$$ Then \be
\lam_1\geq\frac{1}{4}\mu(a)\ctg^2\sqrt{\mu(a)}a\label{5.19}\ee

\end{lem}

{\bf Proof}\,\,\,The Laplace operator comparison theorem yields \be
\Delta_\Pi\rho_\Pi\geq \sqrt{\mu(a)}\ctg\sqrt{\mu(a)}a>0\qfq x\in
M,\,\,\rho_\Pi<a.\label{5.18}\ee

Let $O\subset\subset\Om$ be an open set with a boundary $\pl O$. It
follows from (\ref{5.18})  that \beq\Vol(\pl O)&&=\int_{\pl
O}1d\Ga\geq\int_{\pl O}\<\nabla\rho,\nu\>_\Pi
d\Ga_\Pi=\int_O\Delta_\Pi\rho_\Pi
dg_\Pi\nonumber\\
&&\geq\sqrt{\mu(a)}\ctg\sqrt{\mu(a)}a\Vol(O).\nonumber\eeq Then the
estimate (\ref{5.19}) follows from the Cheeger theorem
(\cite{SYau}). \hfill$\Box$ \\

We have
\begin{thm}\label{t5.2} There is $a>0$ such that, if
$$\Om\subset\{\,x\,|\,x\in
M,\,\,\rho_\Pi(x)<a\,\},$$ then
$$\V_0(\Om)=\{\,0\,\},$$  where $\V_0(\Om)$ is given by
$(\ref{58})$.
\end{thm}

{\bf Proof}\,\,\,For $w\in H^2(\Om)\cap H_0^1(\Om),$ using the
estimate (\ref{5.19}), we have \beq -(w,\,\Delta_\Pi w+\B^*
w)_{L^2_\Pi(\Om)}&&=-(\Delta_\Pi w+\B w,\,
w)_{L^2_\Pi(\Om)}=\||\nabla_\Pi w|\|^2_{L^2_\Pi(\Om)}-(\B
w,w)_{L^2_\Pi(\Om)}\nonumber\\
&&\geq \frac{1}{2}\||\nabla_\Pi
w|\|^2_{L^2_\Pi(\Om)}-C\|w\|^2_{L^2_\Pi(\Om)}\nonumber\\
&&\geq
[\frac{1}{8}\mu(a)\ctg^2\sqrt{\mu(a)}a-C]\|w\|^2_{L^2_\Pi(\Om)}.\nonumber\eeq
The proof is complete.  \hfill$\Box$ \\

{\bf A Elliptic Surface of Revolution}\,\,\,Let $h$ be a smooth
function on $[0,b)$ with $h(0)=0$ such that \be
\frac{1}{s}h''(s)h'(s)>0\qfq s\in[0,b).\label{4.7}\ee Let \be
M=\{\,(x,h(|x|))\,|\,x=(x_1,x_2)\in\R^2,\,\,|x|<b\,\}.\label{4.8*}\ee
We have

\begin{thm}\label{t5.3} Let $o=(0,h(0))$ and let $\Om\subset M$ be a bounded open set which is
 star-shaped with respect to
$o.$ Then \be \V_0(\Om)=\{\,0\,\},\ee where $\V_0(\Om)$ is given
by $(\ref{58})$.
\end{thm}

{\bf Proof}\,\,\,We shall do a careful computation by the formula
(\ref{3.23}). For this end, we make some preparations.

By the formula (\ref{the}), we have
$$\si(\theta)=(\cos\theta,\sin\theta).$$
Let $X,$ $Y$ be vector fields on $\R^2.$ Then
\be\D^2h(X,Y)=[h''(|x|)-\frac{h'(|x|)}{|x|}]\frac{\<X,x\>\<Y,x\>}{|x|^2}+\frac{h'(|x|)}{|x|}\<X,Y\>.
\label{5.26}\ee

We look for a solution to the problem (\ref{the1})-(\ref{the2}) in a
form of $r(t)=\a(t)\si(\theta)$ where $\a(t)>0$ for $t>0$  It is
easy to check from (\ref{5.26}) that such $\a(t)$ is a positive
solution to the problem
\be\left\{\begin{array}{l}\a''(t)+\eta^2(\a(t)\si(\theta))h''(\a(t))h'(\a(t))\a'^2(t)=0\qfq
t>0,\\
\a(0)=0,\quad\a'(0)=1.\end{array}\right.\label{5.27}\ee Moreover,
the solution $\a(t)$ to the problem (\ref{5.27}) is actually the
solution to the problem
\be\left\{\begin{array}{l}\a'(t)=\dfrac{1}{\sqrt{1+h'^2(\a(t))}}\qfq
t>0,\\
\a(0)=0.\end{array}\right.\label{5.28}\ee  Furthermore, a simple
computation shows that $\a(t)$ is also the solution to the problem
(\ref{2.2}). Then
$$\a(t)=f(t),\quad\F(t,\theta)=\Big(r(t), h(r(t))\Big)=\Big(f(t)\si(\theta),h(f(t))\Big)\qfq t\geq0.$$
We obtain \be T=\D_{\pl
t}\F=f'(t)\Big(\si(\theta),h'(f(t))\Big),\quad
E=\frac{1}{f}\D_{\pl\theta}\F=\Big(\dot{\si}(\theta), 0\Big).\ee

We shall prove that the
problem\be\left\{\begin{array}{l}\Delta_\Pi w+\B w=0\qfq x\in\Om,\\
w=0\qfq x\in\Ga,\end{array}\right.\label{5.30}\ee has the unique
zero solution. Let $w$ be a solution to the problem (\ref{5.30}).
By the proof  of Theorem \ref{th3.1}, $(W,w)$ is an infinitesimal
isometry, where $$W=\var T+\phi E,$$
$$\var=-\int_0^tw\Pi_{11}ds,\quad\phi=-w(o)\Pi(\dot{\si}(\theta),\si(\theta))f+\int_0^t\Phi(t,s)P(w)ds.$$
Then \be W(o)=0.\label{531}\ee Let
$$W+wN=(u_1,u_2,u),$$ where $u$ is a solution to the problem
(\ref{3.18}) and $u_1,$ $u_2$ are given by (\ref{3.30}).   Since
$N(o)=(0,0,1)$, it follows from (\ref{531}) that
$$u_1(o)=0,\quad u_2(o)=0,\quad u(o)=w(o).$$
By  (\ref{3.30}), we have
$$z_1=z_2=0.$$ Moreover, using (\ref{3.23}), we obtain  \beq
w/\eta&&=h'(f(t))\cos\theta\int_0^tX_2(t,s)(u)ds-h'(f(t))\sin\theta\int_0^tX_1(t,s)(u)ds\nonumber\\
&&\quad-[1+h'^2(f(t))]u(x)+h'(f(t))\int_0^tu(r(s))h''(f(s))f'(s)ds,\label{5.32}\eeq
where $X_i$ are given by (\ref{3.23*}). On the other, a simple
computation shows that \beq \cos\theta
X_2(t,s)(u)&&=-\{h'(f(s))+[f(s)-f(t)]h''(f(s))]\frac{f'(s)}{f(s)}u_\theta\cos\theta\sin\theta\nonumber\\
&&=\sin\theta X_1(t,s)(u).\label{5.33*} \eeq It follows from
(\ref{5.32}) and (\ref{5.33*}) that \be w\eta=\eta^2
h'(f)\int_0^tu(r(s))h''(f(s))f'(s)ds-u(x)\qfq
x=r(t)\in\Om.\label{5.34}\ee

We now apply the maximum principle to the elliptic problem
(\ref{3.18}) to know that there is $x_0\in\Ga$ such that
$$u(x_0)=\sup_{x\in\Om}u(x).$$ We may assume that $u(x_0)\geq0.$ Otherwise, we consider $-u.$
Let $r(t_0)=x_0.$ Then the formula (\ref{5.34}) yields $$
u(x_0)=\eta^2 h'(f)\int_0^{t_0}u(r(s))h''(f(s))f'(s)ds\leq
u(x_0)\frac{h'^2(f(t_0))}{1+h'^2(f(t_0))},$$  which gives
$u(x_0)=0.$ Next, we consider $-u$ and have $\inf_{x\in\Om}u=0.$
Then $u\equiv0$ on $\Om.$ Finally, we obtain $w\equiv0$ on $\Om$ by
(\ref{5.34}).   \hfill$\Box$ \\

 {\bf A Spherical Cap }\,\,\,Let $M$ be a sphere of constant
curvature $\kappa>0$ with the induced metric $g$ from $\R^3.$ Then
the second fundamental form of $M$ is given by
\be\Pi=\sqrt{\kappa}g.\label{nn4.35}\ee
  Then
$$\sqrt{\kappa}\Delta_\Pi w=\Delta w,\quad\B w=2\sqrt{\kappa}w,$$
where $\Delta$ is the Laplacian of $M$ in the induced metric $g$
from $\R^3.$

 Let
$o\in M$ be given. Let $\rho(x)=\rho(x,o)$ be the distance from
$x\in M$ to $o$ in the induced metric $g$ of $M$.  Set
$$\Om(a)=\{\,x\,|\,x\in M,\,\,\rho(x)<a\,\}\qfq 0<a\leq\frac{\pi}{\sqrt{\kappa}}.$$
Then for $0<a<\dfrac{\pi}{\sqrt{\kappa}},$ $\Om(a)$ is a spherical
cap with a nonempty smooth boundary $$\Gamma(a)=\{\,x\,|\,x\in
M,\,\,\rho(x)=a\,\}.$$ It follows  Theorem \ref{t5.1} that $w\in
H^1_\ib(\Om(a))$ if and only if $w$ satisfies the problem \be \Delta
w+2\kappa w=0\qfq x\in\Om(a).\ee Moreover,
$$\V_0(\Om(a))={\V_0}_*(\Om(a))=\{\,\var\,|\,\Delta\var+2\kappa\var=0,\,\,\var|_{\Ga(a)}=0\,\}.$$

We have

\begin{thm}\label{t4.6}
\be\left\{\begin{array}{l}\V_0(\Om(a))=\{\,0\,\}\qfq 0<a<\dfrac{\pi}{2\sqrt{\kappa}},\\
\V_0(\Om(a))\not=\{\,0\,\}\qfq \dfrac{\pi}{2\sqrt{\kappa}}\leq
a\leq\dfrac{\pi}{\sqrt{\kappa}}.\end{array}\right.\label{5.33}\ee
\end{thm}

{\bf Proof}\,\,\, Let $\lam_1(a)$ be the first positive eigenvalue
of $-\Delta$ on $\Om(a)$ with the Dirichlet boundary condition on
$\Ga(a).$ By  \cite{Ch}, \cite{Ob}, the first positive eigenvalue of
$-\Delta$ of the sphere $M$ without boundary is $2\kappa.$ Since
$C_0^\infty(\Om(a_1))\subset C_0^\infty(\Om(a_2))\subset
C_0^\infty(M)$ for all $0<a_1\leq a_2\leq\frac{\pi}{\sqrt{\kappa}},$
then $H_0^1(\Om(a_1))\subset H^1_0(\Om(a_2))\subset H^1(M)$ in the
following sense: For $h\in H_0^1(\Om(a))$, we define $h=0$ for $x\in
M/\Om(a).$ Then \beq 2\kappa&&=\inf\{\,\frac{\int_M|\nabla
h|^2dg}{\int_Mh^2dg}\,|\,h\in
H^1(M)\,\}\nonumber\\
&&\leq\inf\{\,\frac{\int_{\Om(a_2)}|\nabla h|^2dg}{\int_{\Om(a_2)}
h^2dg}\,|\,h\in
H^1_0(\Om(a_2))\,\}=\lam_1(a_2)\leq\lam_1(a_1).\label{34.3}\eeq

Since $$\Delta\rho=\sqrt{\kappa}\ctg\sqrt{\kappa}\rho\qfq x\in
M,\,\, x\not=o, $$ it is easy to check that the following function
$$\var(x)=\cos\sqrt{\kappa}\rho(x)\qfq x\in M$$ is an eigenfunction of $-\Delta$ of the sphere $M$
without boundary corresponding to the eigenvalue $2\kappa$. Clearly,
$\var$ is also an eigenfunction of $-\Delta$ on
$\Om(\frac{\pi}{2\sqrt{\kappa}})$ with the Dirichlet boundary
condition on $\Ga(\frac{\pi}{2\sqrt{\kappa}})$ corresponding to an
eigenvalue $2\kappa,$ which implies, by (\ref{34.3}), that
$$\lam_1(a)=2\kappa\qfq \frac{\pi}{2\sqrt{\kappa}}\leq
a\leq\frac{\pi}{\sqrt{\kappa}},$$ which means
$$\V_0(\Om(a))\not=\{\,0\,\}  \qfq \frac{\pi}{2\sqrt{\kappa}}\leq
a\leq\frac{\pi}{\sqrt{\kappa}}.$$

Next, we assume that $$0<a<\frac{\pi}{2\sqrt{\kappa}}.$$ Let
$o=(0,0,0)$ and let the semi-sphere
$\Om(\dfrac{\pi}{2\sqrt{\kappa}})$ be given by
$$
\Om(\frac{\pi}{2\sqrt{\kappa}})=\{\,(x,h(|x|))\,|\,x\in\R^2,\,\,|x|<\frac{1}{\sqrt{\kappa}}\,\},$$
where
$$h(s)=\frac{1}{\sqrt{\kappa}}-\sqrt{\frac{1}{\kappa}-s^2}\qfq
s\in[0,\frac{1}{\sqrt{\kappa}}).$$ Since
$$h''(s)h'(s)s^{-1}=\frac{1}{1-\kappa s^2}\qfq s\in[0,\frac{1}{\sqrt{\kappa}}),$$
it follows from Theorem \ref{t5.3} that $\V_0(\Om(a))=\{\,0\,\}$ for
$0<a<\dfrac{\pi}{2\sqrt{\kappa}}.$ \hfill$\Box$

\begin{rem} The relations $(\ref{5.33})$ mean that, for the first eigenvalue $\lam_1(a)$ of $-\Delta$ on
$\Om(a)$ with the Dirichlet boundary condition on $\Ga(a)$,
$$\left\{\begin{array}{l} \lam_1(a)>2\kappa\qfq
0<a<\dfrac{\pi}{2\sqrt{\kappa}},\\
\lam_1(a)=2\kappa\qfq \dfrac{\pi}{2\sqrt{\kappa}}\leq
a\leq\dfrac{\pi}{\sqrt{\kappa}}.\end{array}\right.$$
\end{rem}

\begin{rem}If $a=\dfrac{\pi}{\sqrt{\kappa}},$ then $\Om(a)=M.$ Since a sphere is
rigid, any infinitesimal isometry of $M$ is trivial, see
$\cite{Pog}.$
\end{rem}

\section{Parabolic Surfaces}
\def\theequation{5.\arabic{equation}}
\hskip\parindent A surface $M$ is said to be parabolic if
$$\kappa=0,\quad \Pi\not=0\quad\mbox{for all}\quad x\in M.$$ Let
$M$ be parabolic and orientable. Let $\Om\subset M.$ It follows from
Theorem \ref{th3.1} that $w\in H^1_\is(\Om)$ if and only if $w\in
H^1(\Om)$ solves the problem
$$\<D^2w,Q^*\Pi\>=0\qfq x\in\Om.$$

We assume that there is a vector field $E\in\X(M)$ such that
\be\hat D_EN=0,\quad |E|=0\qfq x\in M.\label{6.1*}\ee Let $p_0\in
M$ be given. We consider a parabolic coordinates $(t,s)$ on $M$ as
follows. Let curves $r$ and $\zeta:$ $\R\rightarrow M$ be given by
$$\left\{\begin{array}{l}\dot{r}(t)=E(r(t))\qfq t\in\R,\\
r(0)=p_0,\end{array}\right.$$ and
$$\left\{\begin{array}{l}\dot{\zeta}(s)=QE(\zeta(s))\qfq s\in\R,\\
\zeta(0)=p_0,\end{array}\right.$$ respectively, where the operator
$Q:\,\,M_p\rightarrow M_p$ for $p\in M$ is given by $(\ref{2.9*}).$
Let two parameters families $\a(t,s)$ and $\b(t,s)$ be given by
$$\left\{\begin{array}{l}\dfrac{\pl\a}{\pl t}(t,s)=E(\a(t,s))\qfq t\in\R,\\
\a(0,s)=\zeta(s),\end{array}\right.$$ and
$$\left\{\begin{array}{l}\dfrac{\pl\b}{\pl s}(t,s)=QE(\b(t,s))\qfq s\in\R,\\
\zeta(t,0)=r(t),\end{array}\right.$$ respectively.  Then
$$\a(t,s)=\b(t,s)\qfq (t,s)\in\R^2,$$
$$\pl t=\frac{\pl\a}{\pl t}(t,s)=E,\quad\pl s=\frac{\pl\b}{\pl
s}(t,s)=Q\frac{\pl\a}{\pl t}(t,s).$$

We have
\begin{thm}\label{t6.1} Let $M$ be a parabolic surface and orientable. Let
$(t,s )$ be the parabolic coordinates on $M$. Then \be
H^1_\ib(M)=\{\,w_0(s)+w_1(s)t\,\,|\,\,w_1,\,\,w_0\in
H^1(\R),\,\,t\in\R\,\,\}.\label{6.2*}\ee
\end{thm}

{\bf Proof}\,\,\,Consider the frame field $E_1=E,$ $E_2=QE.$ By
(\ref{6.1*}), we have \beq D_{\pl t}\pl t&&=\hat D_\pl t\pl
t+\Pi(\pl t,\pl t)N=\frac{\pl^2\a}{\pl
t^2}(t,s)\nonumber\\
&&=\<\frac{\pl^2\a}{\pl t^2}(t,s),\frac{\pl\a}{\pl
t}(t,s)\>\frac{\pl\a}{\pl t}(t,s)+\<\frac{\pl^2\a}{\pl
t^2}(t,s),\frac{\pl\b}{\pl s}(t,s)\>\frac{\pl\b}{\pl
s}(t,s)\nonumber\\
&&\quad+\<\frac{\pl^2\a}{\pl
t^2}(t,s),N\>N\nonumber\\
&&=\frac{1}{2}\frac{\pl}{\pl t}\Big|\frac{\pl\a}{\pl
t}(t,s)\Big|^2\frac{\pl\a}{\pl t}(t,s)+\frac{1}{2}\frac{\pl}{\pl
t}\<\frac{\pl\a}{\pl t}(t,s),Q\frac{\pl\a}{\pl
s}(t,s)\>\frac{\pl\b}{\pl
s}(t,s)\nonumber\\
&&\quad+\frac{\pl}{\pl t}\<\frac{\pl\a}{\pl
t}(t,s),N\>N=0.\label{6.3*}\eeq It follows from (\ref{6.1*}) and
(\ref{6.3*}) that
$$0=\<D^2w,Q^*\Pi\>=D^2w(E,E)\Pi(QE,QE)=\frac{\pl^2w}{\pl
t^2}\Pi(QE,QE).$$ Since $\Pi(QE,QE)\not=0,$ we have the formula
(\ref{6.2*}). \hfill$\Box$\\

{\bf A Cylinder}\,\,\,Let $a>0$ be given. Consider a cylinder
$$M=\{\,(x,z)\,|\,x=(x_1,x_2)\in\R^2,\,\,|x|=a,\,\,z\in\R\,\}.$$
Then $$N=\frac{1}{a}(x,0).$$ Let $E=(0,0,1)=\pl z.$ Then
$$\hat D_EN=0,\quad |E|=1.$$

Consider the parabolic coordinates $(z,\theta)$, given by
$$(x,z)=(a\cos\theta,a\sin\theta,z).$$
 Let $b>0$ be given and let
\be\Om=\{\,(x,z)\,|\,|x|=a,\,\,|z|< b\,\},\quad {\bf
T}=\{\,\,x\,\,|\,\,x\in\R^2,\,\,|x|=a\,\,\}.\label{6.5}\ee Then,
by Theorem \ref{t6.1}, \be
H^1_\is(\Om)=\{\,w_0+w_1z\,|\,w_0,\,w_1\in H^1({\bf
T}),\,\,|z|<b\}.\label{6.1}\ee

\begin{rem} $(i)$ Clearly, $H^1_\is(\Om)\cap C^\infty(\Om)$ is dense in $H^1_\is(\Om).$

$(ii)$ Let $\Om$ be given in $(\ref{6.5})$ with $a=1.$ Let $w_0\in
H^2({\bf T})$ and $w_1\in H^3({\bf T})$ be given. Then an
infinitesimal isometry corresponding to
$w=-w_0'(\theta)+zw_1''(\theta)\in H^1_\is(\Om)$ is given by \beq
V&&=(w_0(\theta)-zw_1'(\theta))\pl\theta+w_1(\theta)\pl
z+wN\nonumber\\
&&=\Big(-w_0(\theta)\sin\theta-w_0'(\theta)\cos\theta,\quad
w_0(\theta)\cos\theta-w_0'(\theta)\sin\theta,\quad
w_1(\theta)\Big)\nonumber\\
&&\quad+z\Big(w_1'(\theta)\sin\theta+w_1''(\theta)\cos\theta,\quad
-w_1'(\theta)\cos\theta+w_1''(\theta)\sin\theta,\quad
0\Big)\nonumber\\
&&=\Big(R(w_0,w_0'),w_1\Big)+z\Big(-R(w_1',w_1''),0\Big)\nonumber\eeq
where
$$R=\left(\begin{array}{cc}-\sin\theta&-\cos\theta\\
\cos\theta&-\sin\theta\end{array}\right).$$\\
\end{rem}

{\bf A Conical Surface}\,\,\,Let $a>0$ and let
$$M=\{\,\,(x,z)\,\,|\,\,|x|=a|z|,\,\,x=(x_1,x_2)\in\R^2,\,\,z\in\R\,\,\}.$$
Then
$$N=\frac{1}{\sqrt{1+a^2}}(\frac{x}{|x|},-1).$$
 Consider the parabolic coordinates $(z,\theta)$, given by
$$(x,z)=z(a\cos\theta,a\sin\theta,1).$$
 Let $b_1,\,\,b_2>0$ be given and let
$$\Om=\{\,\,(x,z)\,\,|\,\,|x|=az,\,\,\,b_1<z< b_2\,\,\},\quad {\bf T}=\{\,\,x\,\,|\,\,x\in\R^2,\,\,|x|=1\,\,\}.$$

Since $\hat D_{\pl z}N=0,$ we have from Theorem \ref{t6.1} \be
H^1_\is(\Om)=\{\,\,w_0+w_1z\,\,|\,\,w_0,\,w_1\in H^1({\bf
T}),\,\,b_1<z<b_2\,\,\}.\label{6.2}\ee

\section{Hyperbolic Surfaces}
\def\theequation{6.\arabic{equation}}
\hskip\parindent A surface $M$ is said to be a hyperbolic surface if
$$\kappa<0\qfq x\in M.$$  Let $M$ be a hyperbolic surface and orientable. We assume
that surface $M$ is given by a family of two parameter curves \be
M=\{\,\,\a(s,\varsigma)\in\R^3\,\,|\,\,(s,\varsigma)\in\R\times{\bf
T}\,\,\},\label{n7.1}\ee which satisfies \be\<\pl s,\pl
\varsigma\>=0,\quad\Pi(\pl s,\pl s)<0,\quad \Pi(\pl
s,\pl\varsigma)=0\qfq x\in M,\label{7.1}\ee where
$${\bf T}=\{\,\,x\in\R^2\,\,|\,\,|x|=1\,\,\}.$$   Let
$\varsigma=(\cos\vartheta,\sin\vartheta)$ for
$\vartheta\in[0,2\pi).$ Then
$$\pl\varsigma=\frac{\pl\a}{\pl\vartheta}(s,\varsigma).$$

We consider the structure of the operator $\<D^2w,Q^*\Pi\>.$ Let
$$E_1=\frac{\pl s}{|\pl s|},\quad
E_2=\frac{\pl\varsigma}{|\pl\varsigma|}.$$ By (\ref{7.1}), we have
\beq |\pl s|^2|\pl\varsigma|^2\<D^2w,Q^*\Pi\>&&=D^2w(\pl s,\pl
s)\Pi(\pl\varsigma,\pl\varsigma)+D^2w(\pl\varsigma,\pl\varsigma)\Pi(\pl
s,\pl s)\nonumber\\
&&=\Pi(\pl\varsigma,\pl\varsigma)w_{ss}+\Pi(\pl
s,\pl s)w_{\vartheta\vartheta}\nonumber\\
&&\quad+\Pi(\pl\varsigma,\pl\varsigma)D_{\pl s}\pl sw+\Pi(\pl s,\pl
s)D_{\pl\varsigma}\pl\varsigma w.\nonumber\eeq

Let $\Om\subset M$ be given by
\be\Om=\{\,\,\a(s,\vartheta)\,\,|\,\,(t,\varsigma)\in(0,b)\times{\bf
T}\,\,\},\label{7.5*}\ee where $b>0$ is given. We fix $o\in\Om$ to
be such that $\Om$ is star-shaped with respect to $o.$

Let
$$a(s,\vartheta)=-\frac{\Pi(\pl
s,\pl s)}{\Pi(\pl\varsigma,\pl\varsigma)}.$$ By Theorem \ref{th3.1},
$w\in H^1_\is(\Om)$ if and only if $w\in H^1(\Om)$ is a solution to
the problem \be
w_{ss}=(a(s,\vartheta)w_\vartheta)_\vartheta+\tilde{\B}w,\label{7.2}\ee
where
\beq\tilde{\B}w&&=-a_\vartheta(s,\vartheta)w_\vartheta-\Pi^{-1}(\pl\varsigma,\pl\varsigma)
[\Pi(\pl\varsigma,\pl\varsigma)D_{\pl s}\pl sw+\Pi(\pl s,\pl
s)D_{\pl\varsigma}\pl\varsigma w]\nonumber\\
&&\quad-|\pl
s|^{-2}|\pl\varsigma|^{-2}\Pi^{-1}(\pl\varsigma,\pl\varsigma)[w\kappa\tr\Pi+\kappa_1\int_0^tw\Pi_{11}ds
+\kappa_2\int_0^t\Phi(t,s)P(w)(s)ds\nonumber\\
&&\quad+w(o)\Pi(\si(\theta),\dot{\si}(\theta))\kappa_2f],\label{7.3}\eeq
where   the four factor of the third term in the right hand side of
(\ref{7.3}) is given by Theorem \ref{th3.1}. Clearly, the linear
operator\be \tilde{\B}:\quad H^1(\Om)\rightarrow
L^2(\Om),\label{7.4}\ee is bounded.

 We introduce a family
of self-adjoint operators on $L^2({\bf T})$ by
$$A(s)u=-\Big(a(s,\vartheta)u_\vartheta\Big)_\vartheta,\quad D(A(s))=H^2({\bf T})\qfq s\in[0,b].$$
Consider a family of operators  on $H^1({\bf T})\times L^2({\bf T})$
$${\bf A}(s)=\left(\begin{array}{cc}0&I\\
A(s)&0\end{array}\right),\quad D(\tilde{A}(s))=H^2({\bf T})\times
H^1({\bf T})\qfq s\in[0,b].$$ Let
$$\grave{H}^m({\bf T})=\{\,\,u\in H^m({\bf
T})\,\,|\,\,(u,1)_{L^2({\bf T})}=0\,\,\}\qfq m=0,\,\,1.$$

\begin{lem}\label{l7.1} The operator family $\{\,\,{\bf A}(s)\,\,\}_{0\leq s\leq b}$ generates a unique
evolution system ${\bf U}(s,\lam)$ for $0\leq\lam\leq s\leq b$ on
$H^1({\bf T})\times \grave{L}^2({\bf T}).$ In particular, there
exist constants $C(b)>0$ and $\omega(b)>0$ such that \be\|{\bf
U}(s,\lam)\|\leq C(b)e^{\omega(b)(s-\lam)}\qfq 0\leq\lam\leq s\leq
b.\label{7.9}\ee
\end{lem}

{\bf Proof}\,\,\,We introduce the equivalent norms on $H^m({\bf
T})\times H^{m-1}({\bf T})$ by
$$\|(u,v)\|_m=\Big(\frac{1}{2\pi}(u,1)^2+\sup_{0\leq s\leq b}[(A^m(s)u,u)_{L^2({\bf T})}
+(A^{m-1}(s)v,v)_{L^2({\bf T})}]\Big)^{1/2}$$ for $(u,v)\in
H^m({\bf T})\times H^{m-1}({\bf T}),$ $m=1,$ $2,$ and $3$,
respectively, where $H^0({\bf T})=L^2({\bf T})$ and $A^0(s)=I.$
Note that $\Big(H^m({\bf T})\times H^{m-1}({\bf
T}),\quad\|\cdot\|_m\Big)$ is not a Hilbert space in general.

Let $$0<\lam_1(s)\leq\lam_2(s)\leq\cdots\leq\lam_k(s)\leq\cdots$$
be all positive eigenvalues of $A(s)$ and let their corresponding
eigenfunctions be $\{\,\var_k\,\}$  such that
$\{(2\pi)^{-1/2},\,\,\var_k\,\}$ forms an orthonormal basis of
$L^2({\bf T}).$  Then for each $s\in[0,b]$ the operator ${\bf
A}(s)$ generates a $C_0$ group semigroup ${\bf S}_s(t)$ on
$H^1({\bf T})\times L^2({\bf T})$, given by
\be {\bf S}_s(t)\left(\begin{array}{c}u\\v\end{array}\right)=\left(\begin{array}{c}c_0\\
0\end{array}\right)+
\sum_k a_ke^{\sqrt{\lam_k}it}\left(\begin{array}{c}\frac{1}{\sqrt{\lam_k}}\var_k\\
i\var_k\end{array}\right)+\sum_k b_ke^{-\sqrt{\lam_k}it}\left(\begin{array}{c}\frac{1}{\sqrt{\lam_k}}\var_k\\
-i\var_k\end{array}\right),\label{7.5} \ee where
$$c_0=\frac{1}{2\pi}(u,1)_{L^2},\quad \left\{\begin{array}{l}a_k+b_k=\sqrt{\lam_k}(u,\var_k)_{L^2},\\
a_k-b_k=-i(v,\var_k)_{L^2},\end{array}\right.$$ and $ (u,v)\in
H^1({\bf T})\times \grave{L}^2({\bf T})$ is real.

Let ${\bf S}_s(t)(u,v)=\Big({\bf S}_{s1}(t),{\bf S}_{s2}(t)\Big).$
It follows from (\ref{7.5}) that \beq &&\Big(A(s){\bf
S}_{s1}(t),{\bf
S}_{s1}(t)\Big)_{L^2}+\|{\bf S}_{s2}(t)\|^2_{L^2}=\sum_k(|a_ke^{\sqrt{\lam_k}it}+b_ke^{-\sqrt{\lam_k}it}|^2
+|a_ke^{\sqrt{\lam_k}it}-b_ke^{-\sqrt{\lam_k}it}|^2)\nonumber\\
&&=2\sum_k(|a_k|^2+|b_k|^2)=\frac{1}{2}\sum_k[\lam_k(u,\var_k)_{L^2}^2+(v,\var_k)_{L^2}^2]
=\frac{1}{2}[(A(s)u,u)_{L^2}+(v,v)_{L^2}],\nonumber\eeq which
yield
$$\|{\bf S}_s(t)(u,v)\|_1\leq \|(u,v)\|_1\qfq (u,v)\in H^1({\bf T})\times \grave{L}^2({\bf T}),\quad s\in[0,b].$$
Similarly, we have
$$\|{\bf S}_s(t)(u,v)\|_2\leq \|(u,v)\|_2\qfq (u,v)\in H^2({\bf T})\times \grave{H}^1({\bf T}),\quad s\in[0,b].$$

Then the proof is complete by Theorem 3.1 of Chapter 5 in
\cite{Pazy}. \hfill$\Box$

\begin{thm} \label{t7.1} For any $w_0\in H^1({\bf T})$ and $w_1\in
L^2({\bf T})$ with $(w_1,1)_{L^2({\bf T})}=0,$ there is a unique
$w\in H^1_\is(\Om)$ such that \be
w(0,\vartheta)=w_0(\vartheta),\quad
w_s(0,\vartheta)=w_1(\vartheta).\label{7.6}\ee
\end{thm}

{\bf Proof}\,\,\, The problem (\ref{7.2}) and (\ref{7.6}) is
equivalent to the first order system\be
\left\{\begin{array}{l}\dfrac{\pl}{\pl
s}\left(\begin{array}{c}w\\ w_s\end{array}\right)=\left(\begin{array}{cc}0&I\\
A(s)&0\end{array}\right)\left(\begin{array}{c}w\\
w_s\end{array}\right)+\left(\begin{array}{c}0\\\tilde{\B}w\end{array}\right)\qfq (s,\varsigma)\in(0,b)\times{\bf T},\\
(w(0),w_s(0))=(w_0,w_1)\qfq \varsigma\in{\bf
T}.\end{array}\right.\label{7.7}\ee

Let $(w_0,w_1)\in H^1({\bf T})\times \grave{L}^2({\bf T})$ be
given. Let $b\geq\eta>0$ be given. Consider a Banach space, given
by
$${\bf X}_\eta= L^1\Big(0,\eta;H^1((0,\eta)\times{\bf T})\Big).$$
Consider a linear operator ${\bf F}:$ ${\bf X}_\eta\rightarrow{\bf
X}_\eta$, given as follows. For $u\in{\bf X}_\eta$, ${\bf F}u$ is
defined as the first component of \be{\bf
U}(s,0)(w_0,w_1)+\int_0^s{\bf U}(s,\lam)(0,\tilde{\B}u)d\lam\qfq
0\leq s\leq\eta,\label{7.8}\ee where the operator $\tilde{\B}$ is
given by (\ref{7.3}) and ${\bf U}$ is the evolution system in
Lemma \ref{l7.1}. Then a solution $w\in {\bf X}_\eta$ to the
problem (\ref{7.7}) for $0\leq s\leq\eta$ if and only if $w$ is a
fixed point of the operator ${\bf F}$ in ${\bf X}_\eta.$

 Let $u_0$ be the first component of the
first term in (\ref{7.8}). For $u\in {\bf X}\eta$, let ${\bf G}u$
be the first component of the second term in (\ref{7.8}). Then
$${\bf F}u=u_0+{\bf G}u\qfq u\in{\bf X}_\eta.$$
Moreover, we have, by (\ref{7.9}) and (\ref{7.4}), \beq&&\|({\bf
G}u)(s)\|_{H^1((0,\eta)\times{\bf T})}=\|\int_0^s{\bf
U}(s,\lam)(0,\tilde{\B}u)d\lam\|_{H^1({\bf T})\times
L^2({\bf T})}\nonumber\\
&&\leq C(b)\int_0^se^{\omega(b)(s-\lam)}d\lam\sup_{0\leq
s\leq\eta}\|\tilde{\B}u\|_{L^2({\bf T})}\nonumber\\
&&\leq \frac{C(b)}{\omega(b)}(e^{\omega(b)\eta}-1)\|u\|_{{\bf
X}_\eta}\qfq 0\leq s\leq\eta, \quad u\in {\bf
X}_\eta.\nonumber\eeq Then for $\eta>0$ small, ${\bf G}$ is a
strictly contractive map on ${\bf X}_\eta$ which implies that
${\bf F}$ has a unique fixed point $w\in{\bf X}_\eta,$ given by
$$w=\sum_{k=0}^\infty {\bf G}^ku_0\qfq 0\leq s\leq\eta.$$
 Then $w$ is a
solution to the problem (\ref{7.7}) for $0\leq s\leq\eta.$ Moreover,
the solution $w$ can be extended to $s\in[0,b]$ since the problem
(\ref{7.7}) is linear. \hfill$\Box$

\begin{rem}  Let $\Om$ be given by $(\ref{7.5*})$ and let $o\in\Om$
be fixed. In general
$$\Om\cap\exp_oC(o)\not=\emptyset,$$ where $\exp_oC(o)$ is the cut
locus of $o.$ Then the operator $\t{\B}$, given by $(\ref{7.3}),$
does not map $C^\infty(\Om)$ into $C^\infty(\Om).$ Then Theorem
$\ref{t7.1}$ does not imply density results.
\end{rem}

Let $M$ be given by (\ref{n7.1}). To obtain density results, we need
to choose $\Om$ such that $\Om\subset\exp_o\Sigma(o).$ For this end,
we let \be \Om=\{\,\,\a(s,\vartheta)\in
M\,\,|\,\,s\in(0,b),\,\,\vartheta\in(0,\vartheta_0)\,\,\},\label{n6.12}\ee
where $b>0$ and $0<\vartheta_0<2\pi.$ Let $o\in\Om$ be fixed. Let
$$M_1=M/\{\,\a(s,\vartheta_1)\,\,|\,s\in\R\,\,\},\quad \vartheta_0<\vartheta_1<2\pi.$$   Since
$(M_1,g)$ is simply connected and curvature negative, we have
$$\exp_o\Sigma(o)=M_1.$$  In this sense
\be\bar{\Om}\subset\exp_o\Sigma(o).\label{n6.13}\ee

By similar arguments as for Theorem \ref{t7.1}, we have

\begin{thm} Let $\Om$ be given by $(\ref{n6.12}).$ For $h_1,\,\,h_2\in
H^1(0,b),$ $w_0\in H^1(0,\vartheta_0),$ and $w_1\in
L^2(0,\vartheta_0)$ given, there is a unique $w\in H^1_\is(\Om)$
such that
$$w(s,0)=h_1(s),\quad w(s,\vartheta_0)=h_2(s),\quad
w(0,\vartheta)=w_0(\vartheta),\quad
w_s(0,\vartheta)=w_1(\vartheta).$$
\end{thm}

By similar arguments as in Theorem \ref{tn4.3}, it follows from
Theorem \ref{7.1} and the relation (\ref{n6.13}) that

\begin{thm}\label{nt6.3} Let $\Om$ be given by
$(\ref{n6.12}).$ Then the strong $H^1(\Om)$ closure of
$$H^1_{\ib}(\Om)\cap C^\infty(\Om)$$ agrees with $H^1_{\ib}(\Om).$
\end{thm}

We present two examples which satisfy the assumptions (\ref{7.1}) to
end this section.\\

 {\bf A Segment Surface of Revolution}\,\,\, Let
$$
\t{M}=\{\,\,\a(r,\vartheta)\,\,|\,\,r\geq0,\,\,\vartheta\in(0,2\pi]\,\,\},$$
where
$$\a(r,\vartheta)=(r\cos\vartheta,r\sin\vartheta,\log(1+r^2)).$$Then
$$\pl r=(\cos\vartheta,\sin\vartheta,\frac{2r}{1+r^2}),\quad\pl\vartheta=r(-\sin\vartheta,\cos\vartheta,0),$$
$$N=\frac{2r}{\sqrt{1+6r^2+r^4}}(\cos\vartheta,\sin\vartheta,-\frac{1+r^2}{2r}).$$

We have \beq \Pi(\pl r,\pl
r)&&=\frac{2(1-r^2)}{(1+r^2)\sqrt{1+6r^2+r^4}},\nonumber\eeq
$$\Pi(\pl\vartheta,\pl\vartheta)=\frac{2r^2}{\sqrt{1+6r^2+r^4}},\quad \Pi(\pl r,\pl\vartheta)=0,$$
$$\kappa=\frac{1}{|\pl r|^2|\pl\vartheta|^2}\Pi(\pl r,\pl r)\Pi(\pl\vartheta,\pl\vartheta)=
\frac{4(1-r^4)}{(1+6r^2+r^4)^2}\left\{\begin{array}{l}>0\qfq
0<r<1;\\
=0\qfq r=1;\\
<0\qfq r>1.\end{array}\right.$$

We let
$$M=\{\,\,\a(r,\vartheta)\,\,|\,\,r>1,\,\,\vartheta\in[0,2\pi)\,\,\}.$$
Then the assumptions (\ref{7.1}) hold true. \\

{\bf A Hyperboloid of one Sheet}\,\,\, Let
$$M=\{\,\,(x,z)\,\,|\,\,x=(x_1,x_2)\in\R^2,\,\,z^2+1=x_1^2+x_2^2\,\,\}.$$
Consider a family of two parameter curves
$$\a(r,\vartheta)=\Big(r\cos\vartheta,r\sin\vartheta, \sqrt{r^2-1}\Big)\qfq r>1,\quad\vartheta\in[0,2\pi).$$
Then
$$\pl
r=(\cos\vartheta,\sin\vartheta,\frac{r}{\sqrt{r^2-1}}),\quad\pl\vartheta=r(-\sin\vartheta,\cos\vartheta,0),
$$
$$ N=\eta(\cos\vartheta,\sin\vartheta,-\frac{\sqrt{r^2-1}}{r}),$$
$$\eta=\frac{r}{\sqrt{2r^2-1}},$$
$$\Pi(\pl r,\pl r)=-\frac{\eta}{r(r^2-1)},\quad \Pi(\pl r,\pl\vartheta)=0,\quad
\Pi(\pl\vartheta,\pl\vartheta)=r\eta.$$ The assumptions (\ref{7.1})
hold.

\section{Bending of Shells}
\def\theequation{7.\arabic{equation}}
\hskip\parindent We shall apply the theories in Sections 3-6 to
the limit energy functionals of the $\Ga$-convergence to reduce
bending of shells to a one-dimensional problem in the elliptic
case, or parabolic case, or hyperbolic case, respectively.

 Let $M$ be a connected, oriented surface in $\R^3$
with the normal field $N$. Suppose that $g$ is the induced metric
of the surface $M$ from the standard metric of $\R^3$.    A family
$\{\S^h\}_{h>0}$ of shells of small thickness $h$ around $\Om$ is
given through
$$\S^h =\{ \,p\,|\, p=x+zN(x),\,\, x\in\Omega,\,\,-h/2 < z < h/2\,\},
\quad 0<h<h_0.$$ The projection onto $\Om$ along $N$ will be
denoted by $\pi.$ We will assume that $0<h < h_0$, with $h_0$
sufficiently small to have $\pi$ well defined on each $\S^h.$

To a deformation $u\in W^{1,2}(\S^h,\R^3)$, we associate its elastic
energy (scaled per unit thickness): \be  E^h(u) = \frac{1}{h}
\int_{\S^h} W(\hat{\nabla}u)dp, \label{8.1}\ee where $\hat{\nabla}$
denotes the gradient of the Euclidean space $\R^3.$ Here, the stored
energy density $W:$ $\R^{3\times3}\rightarrow [0,\infty]$ is assumed
to be $C^2$ in a neighborhood of $\SO(3),$ and to satisfy the
following normalization, frame indifference and nondegeneracy
conditions
$$\forall\,F\in\R^{3\times3},\,\,\forall\,R\in\SO(3),\quad W(R)=0,\quad
W(RF)=W(F),$$
$$W(F)\geq C\dist^2(F,\SO(3))$$
(with a uniform constant $C>0$). In the study of the elastic
properties of thin shells $\S^h,$ a crucial step is to describe the
limiting behavior, as $h\rightarrow0,$ of minimizers $u^h$ to the
total energy functional \be
J(u)=E^h(u)-\frac{1}{h}\int_{\S^h}\<f^h,u\>dp,\label{8.2}\ee subject
to applied forces $f^h.$ It can be shown that if the forces $f^h$
scale like $h^\a,$ then $E^h(u) \sim h^\b$ where $\b=\a$ if
$0\leq\a\leq2$ and $\b=2\a-2$ if $\a>2.$ The main part of the
analysis consists, therefore, of characterizing the limiting
behavior of the scaled energy functionals $\dfrac{E^h}{h^\b},$ or
more generally, that of $\dfrac{E^h}{e^h},$ where $e^h$ is a given
sequence of positive numbers obeying a prescribed scaling law.

The first result in this framework is due to \cite{LeRa1}, who
studied the scaling $\b=0.$ This leads to a membrane shell model
with energy depending only on stretching and shearing of the
mid-surface. The case $\b=2$ has been analyzed in \cite{FrMoMu} and
it corresponds to geometrically nonlinear bending theory, where the
only admissible deformations are the isometries of the mid-surface,
while the energy expresses the total change of curvature produced by
the deformation.

In \cite{LeMoPa}, the limiting model has been identified for the
range of scalings $\b\geq4$, based on some estimates in
\cite{FrMoMu1}. In these cases, the admissible deformations u are
only those which are close to a rigid motion R and whose first order
term in the expansion of $u-R$ with respect to $h$ is given by $RV$,
where $V\in \IS^1(\Om,\R^3)$ is an infinitesimal isometry on $\Om.$

Let $V\in\IS^1(\Om,\R^3).$ Then there exists a matrix $A$ such that
\be A^\tau(x)=-A(x),\quad \hat{D}_XV=A(x)X,\quad X\in M_x,\quad
x\in\Om.\label{8.8}\ee

For $\b>4$ the limiting energy is given only by a bending term, that
is, the first order change in the second fundamental form of $\Om,$
produced by $V$, \be
I(V)=\frac{1}{24}\int_\Om\Q_2\Big(x,\,\,\Xi(V)\Big)dg \qfq
V\in\IS^1(\Om,\R^3),\label{8.3}\ee where
\be\Xi(V)=(\hat{D}^*(AN)-A\Pi)_{\tan},\label{8.5}\ee and corresponds
to the linear pure bending theory derived in \cite{Ci} from
linearized elasticity. In (\ref{8.5}), $\hat D^*(AN)$ is the
transpose of $\hat D(AN)$, given by
$$\hat D^*(AN)(\tau,\eta)=\<\hat D_{\tau}(AN),\eta\>\qfq
\tau,\,\,\eta\in M_x,\,\,x\in\Om.$$

In (\ref{8.3}), the quadratic forms $\Q_2(x, \cdot)$ are defined as
follows: $$\Q_2(x, F_{\tan})=\min_{a\in\R^3}\Q_3(F+a\otimes N),\quad
\Q_3(F) = D^2W(I)(F, F).$$ The form $\Q_3$ is defined for all
$F\in\R^{3\times3},$ while $\Q_2(x, \cdot)$ for a given $x\in\Om$ is
defined on tangential minors
$F_{\tan}=(\<F\tau,\eta\>)_{\tau,\,\,\eta\in M_x}$ of such matrices.

For $\b=4$ the $\Ga$-limit, which turns out to be the generalization
of the von K¨¢rm¨¢n functional \cite{FrMoMu2} to shells, also
contains a stretching term measuring the second order change in the
metric of $\Om$,
$$\t{I}(V,B_{\tan})=\frac{1}{2}\int_\Om\Q_2\Big(x,\,\,B_{\tan}-\frac{1}{2}(A^2)_{\tan}\Big)dg
+\frac{1}{24}\int_\Om\Q_2\Big(x,\,\,\Xi(V)\Big)dg$$ for
$V\in\IS^1(\Om,\R^3).$ This involves a symmetric matrix field
$B_{\tan}$ belonging to the finite strain space
$${\bf B}=\Big\{\,\,L^2-\lim_{h\rightarrow0}\sym\hat{\nabla}w^h\,\,|\,\,w\in W^{1,2}(\Om,\R^3)\,\,\Big\},$$
where
$$\sym\t{\nabla}w(\tau,\eta)=\frac{1}{2}(\<\t{\nabla}w\tau,\eta\>+\<\t{\nabla}w\eta,\tau\>)
\qfq \tau,\,\,\eta\in M_x,\quad x\in\Om.$$   The space ${\bf B}$
emerges as well in the context of linear elasticity and
ill-inhibited surfaces \cite{GeSa, SaPa}.

It was further shown in \cite{LeMoPa} that for a certain class of
surfaces, referred to as approximately robust surfaces, the limiting
energy for $\b=4$ reduces to the purely linear bending functional
(\ref{8.3}). Elliptic  surfaces happen to belong to this class
\cite{LeMoPa}.

Moreover, \cite{LeMoPa2} has proved that   the limit energy of the
range of scalings $2<\b<4$ for elliptic surfaces is still given by
(\ref{8.3}).

Here we focus on the limit energy (\ref{8.3}) and reduce it from
over the space $\IS^1(\Om,\R^3)$ to over the space $H^1_\is(\Om)$
to give mathematical formulas, as in \cite{Yao}.

Let $T_0\in T^2(M)$ be the third fundamental form of surface $M$,
given by
$$ T_0(\tau,\eta)=\<\hat D_\tau N,\hat D_\eta N\>\qfq
\tau,\,\,\eta\in M_x,\quad x\in M.$$

We now describe the limiting energy formula (\ref{8.3}) in the
common denotation in Riemannian geometry. For simplicity, we
restrict ourselves  to the case when the stored-energy function is
isotropic (that is to say, $W(F)=W(R_1FR_2)$ for all $F\in
M^{3\times3}$ and all $R_1,$ $R_2\in\SO(3)$). In this case, the
second derivative of $W$ at the identity is
$$D^2W(I)(A,A)=2\mu|E|^2+\lam(\hat\tr E)^2,\quad E=\frac{A+A^\tau}{2},$$
for some constants $\mu,$ $\lam\in\R.$

\begin{lem}
 Let $\mu>0$ and $2\mu+\lam>0.$ For $G\in T^2(\Om)$ symmetric,
\be\Q_2(x,G)=2\mu|G|_{T^2_x}^2+\frac{\lam\mu}{\mu+\lam/2}\tr^2G\qfq
x\in\Om,\label{8.7}\ee where $|\cdot|_{T_x^2}$ is given by
$(\ref{8.4})$ and $\tr$ is the trace in the induced metric $g.$
\end{lem}

{\bf Proof}\,\,\,Let $x\in \Om$ be given and let
$${\bf F}=\{\,\,F\,\,|\,\,F\in
M^{3\times3},\,\,\mbox{symmetric}\,\,\},\quad {\bf
F}_0=\{\,\,a\otimes N+N\otimes a\,\,|\,\,a\in\R^3\,\,\}.$$ We
introduce an inner product on ${\bf F}$ by
$$\<F_1,F_2\>_*=2\mu\<F_1,F_2\>+\lam\hat\tr F_1\hat\tr F_2\qfq
F_1,\,\,F_2\in {\bf F}.$$ Then $({\bf F},\<\cdot,\cdot\>_*)$ is an
inner product space.

Let $e_1,$ $e_2$ be an orthonormal basis of $M_x.$ Then $F_1,$
$F_2,$ $F_3$ forms an orthonormal basis of $({\bf
F}_0,\<\cdot,\cdot\>_*),$ where
$$F_i=\frac{e_i\otimes N+N\otimes e_i}{2\sqrt{\mu}},\quad F_3=\frac{1}{\sqrt{2\mu+\lam}}N\otimes N.$$
 Then for $F\in M^{3\times3}$ symmetric with $G=F_{\tan},$ we have \beq
\Q_2(x,G)&&=\min_{a\in\R^3}\Q_3(F+a\otimes
N)=|F-\sum_{i=1}^3\<F,\quad
F_i\>_*F_i|^2_*\nonumber\\
&&=2\mu|F-\sum_{i=1}^3\<F,\quad F_i\>_*F_i|^2+\lam(\hat \tr
F-\<F,F_3\>_*\hat\tr F_3)^2\nonumber\\
&&=2\mu[\sum_{ij=1}^2\<F,e_i\otimes
e_j\>^2+(\frac{\lam}{2\mu+\lam})^2(\sum_{i=1}^2\<F,e_i\otimes
e_i\>)^2]\nonumber\\
&&\quad+\lam(\frac{2\mu}{2\mu+\lam})^2(\sum_{i=1}^2\<F,e_i\otimes
e_i\>)^2\nonumber\\
&&=2\mu|G|_{T^2_x}^2+\frac{2\mu\lam}{2\mu+\lam}\tr^2G.\nonumber\eeq
\hfill$\Box$

Let $k$ be a nonnegative integer and let $T\in T^k(\Om)$ be a
$k$th-order tensor field on $\Om.$ The internal product of $X$
with $T$ is a $k-1$-th order tensor field $\i(X)T,$ defined by
\be\i(X)T(X_1,\cdots,X_{k-1})=T(X,X_1,\cdots,X_{k-1})\qfq
X_1,\,\,\cdots,\,\,X_{k-1}\in\X(M).\label{8.9}\ee

\begin{lem}
Let $V\in\IS^1(\Om,\R^3)$ with $V=W+wN.$ Then
\be\Xi(V)=\i(W)D\Pi+\Pi(D_\cdot W,\cdot)+\Pi(\cdot,D_\cdot
W)+wT_0-D^2w,\label{8.6}\ee  where $\Xi(V)$ is given by
$(\ref{8.5})$, $D$ is the Levi-Civita connection of the induced
metric $g$, and $\cdot$ denotes the position of variables.
\end{lem}

{\bf Proof}\,\,\,It follows from (\ref{8.8}) that
$$AX=D_XW+w\hat D_XN+[X(w)-\Pi(W,X)]N\qfq X\in M_x,\,\,x\in\Om.$$
Then \be\<AN,X\>=-\<N,AX\>=\Pi(W,X)-W(w)\qfq X\in
M_x,\,\,x\in\Om.\label{8.10}\ee Since $\<AN,N\>=0,$ the identity
(\ref{8.10}) yields
$$AN=\i(W)\Pi-Dw.$$

Let $x\in\Om$ be given. We compute the identity (\ref{8.6}) at the
point $x$. Let $e_1,$ $e_2$ be an orthonormal basis of $M_x$ such
that
$$\hat D_{e_i}N=\lam_iN,\quad \lam_i=\Pi(e_i,e_i)\qfq i=1,\,\,2.$$
Let $E_1,$ $E_2$ be a frame field normal at $x$ such that
$$E_i=e_i\quad\mbox{at}\quad x\qfq i=1,\,\,2.$$
We have at $x$ \beq \Xi(V)(\tau_i,\tau_j)&&=\hat
D^*(AN)(\tau_i,\tau_j)-A\Pi(\tau_i,\tau_j)=\<D_{\tau_i}(AN),\tau_j\>-\<A\Pi\tau_i,\tau_j\>\nonumber\\
&&=D\Pi(W,\tau_j,\tau_i)+\Pi(D_{\tau_i}W,\tau_j)-\tau_i\tau_j(w)+\lam_i\<\tau_i,A\tau_j\>\nonumber\\
&&=D\Pi(W,\tau_j,\tau_i)+\Pi(D_{\tau_i}W,\tau_j)+\Pi(D_{\tau_j}W,\tau_i)+w\lam_i\lam_j-\tau_i\tau_j(w),\nonumber\eeq
which yields the identity (\ref{8.6}). \hfill$\Box$

\begin{rem} The identity $(\ref{8.6})$ shows that the tensor field $\Xi(V)$, gievn
by $(\ref{8.3}),$ is exactly the  change of the linearized curvature
tensor of the middle surface $\Om,$ introduced by $\cite{Ko},$ also
see $\cite{BeBo, CY}$ or $\cite{Yao},$  in the case of infinitesimal
deformations.

Consider a deformation $\var:$ $\Om\rightarrow\R^3.$ After the
deformation, the middle surface becomes
$$\bar{\Om}=\{\,\,\var(x)\,\,|\,\,x\in\Om\,\,\}.$$ Let $\bar{\Pi}$
be the second fundamental form of $\bar{\Om}.$ Then the change of
curvature tensor of the middle surface is defined by
$$G=\var^*\bar{\Pi}-\Pi,$$ which is a $2$-th tensor field on $\Om.$

Consider a small deformation $$\var(x)=x+V(x)\qfq x\in\Om$$ with
$V=W+wN\in\IS^1(\Om,\R^3).$ Let $\bar{N}$ be the normal of
$\bar{\Om}.$ After linearization $(\cite{Yao})$, we have
$$\bar{N}(\var(x))=\i(W)\Pi-Dw+N\qfq x\in\Om.$$

Let $x\in\Om$ be given and let $E_1,$ $E_2$ be a frame field normal
at $x$ with the positive orientation. Then \beq
&&\mbox{Linearization}\,(\var^*\bar{\Pi}-\Pi)(E_i,E_j)=\mbox{Linearization}\,[\<\hat
D_{\var_*E_i}\bar{N},\var_*E_j\>-\Pi(E_i,E_j)]\nonumber\\
&&=\<\hat D_{E_i}[\i(W)\Pi-Dw],\quad E_j\>+\<\hat
D_{E_i}N,\quad\var_*E_j\>-\Pi(E_i,E_j)\nonumber\\
&&=D\Pi(W,E_i,E_j)+\Pi(D_{E_i}W,E_j)+\Pi(D_{E_j}W,E_i)\nonumber\\
&&\quad+wT_0(E_i,E_j)-D^2w(E_i,E_j)\quad\mbox{at}\quad
x,\nonumber\eeq  that is, by $(\ref{8.6}),$
$$\Xi(V)=\mbox{Linearization}\,(\var^*\bar{\Pi}-\Pi).$$
\end{rem}

Let $\Om\subset M$ be elliptic and star-shaped  with respect to
$o\in\Om.$ We further assume that for any $\psi\in H^{1/2}(\Ga)$
the problem (\ref{5.9}) has a unique solution $w=\lambda(\psi)\in
H^1(\Om).$ By Theorem \ref{t5.1}, there is a unique
$W=\Lambda(\psi)\in H^1(\Om,T)$ which is perpendicular to
$H^1_\kf(\Om,T)$ such that $V=\Lambda(\psi)+\lam(\psi)N$ is an
infinitesimal isometry. Then for any $V\in\IS^1(\Om,\R^3)$, we
have a formula in the form of
$$V=W+\Lambda(\psi)+\lam(\psi)N\qfq W\in H^1_\kf(\Om,T),\,\,\psi\in
H^{1/2}(\Ga).$$ By Theorem \ref{tn3.1}, $\dim
H^1_\kf(\Om,T)\leq3.$ Then the limit energy (\ref{8.3}) of the
$\Ga$-convergence becomes a functional over  a one-dimensional
space \be I(V)=\t{I}(\a,\psi)\qfq (\a,\psi)\in\R^3\times
H^1(\Ga).\label{8.10}\ee Similar situations happen when the meddle
surface $\Om$ is parabolic or hyperbolic. It follows from Theorems
\ref{t5.1},  \ref{t6.1}, and \ref{t7.1} that
\begin{thm}\label{t8.1}
 Let the meddle surface $\Om$ be elliptic, or parabolic, or
 hyperbolic. Then the limit energy formula $(\ref{8.3})$ of the $\Ga$-convergence
 reduces to be a one-dimensional problem.
\end{thm}

We shall write out explicit formulas of (\ref{8.10}) for spherical
shells and cylinder shells, respectively, before ending this
section.\\

{\bf Bending of Spherical Shells}\,\,\, Let $M$ be the sphere of
curvature $\kappa>0$ and let $g$ be the induced metric of $M$ from
$\R^3.$ Then  the third fundamental form of $M$ is given by
$$T_0=\kappa g.$$

Let $o\in M$ be fixed. Let $\rho(x)=\rho(x,o)$ be the distance
function from $x\in M$ to $o$ in the induced metric $g.$ For
$0<a\leq\dfrac{\pi}{\sqrt{\kappa}},$ let
\be\Om(a)=\{\,\,x\,\,|\,\,x\in
M,\,\,\rho(x)<a\,\,\},\quad\Ga(a)=\{\,\,x\,\,|\,\,x\in
M,\,\,\rho(x)=a\,\,\}.\label{8.17}\ee

Let $V=W+wN$ be an infinitesimal isometry on $\Om(a).$ By the
formulas (\ref{8.6}) and (\ref{nn4.35}), we have \beq
\Xi(V)&&=\sqrt{\kappa}(DW+D^*W)+\kappa wg-D^2w=-\kappa
wg-D^2w.\label{8.11}\eeq In particular, for $V=W\in H^1_\kf(\Om,T)$
a Killing field,
$$\Xi(V)=0.$$

By Theorem \ref{t5.1}, $w\in H^1_\is(\Om)$ if and only if $w$ solves
the problem \be\left\{\begin{array}{l}\Delta w+2\kappa w=0\qfq x\in\Om(a),\\
w=\psi\qfq x\in\Ga(a).\end{array}\right.\label{8.12}\ee Then it
follows from (\ref{8.11}) and (\ref{8.12}) that \be \tr\Xi(V)=0\qfq
x\in\Om.\label{8.13}\ee

Furthermore, we have

\begin{lem} Let $V=W+wN$ be an infinitesimal isometry with $w\in
H^1_\is(\Om).$ Then \be
|\Xi(V)|_{T^2_x}^2=\frac{1}{2}\Delta|Dw|^2+\kappa\div w\nabla w \qfq
x\in\Om.\label{8.14}\ee
\end{lem}

{\bf Proof}\,\,\,Recall that  the Weitzenb$\ddot{o}$ck
formula(Theorem 1.27 in \cite{Yao}) reads \be
|D^2w|_{T^2_x}^2=\frac{1}{2}\Delta|Dw|^2+\<{\bf
\Delta}Dw,Dw\>-\Ric(Dw,Dw)\qfq x\in\Om,\label{8.15}\ee where ${\bf
\Delta}$ is the Hodge-Laplacian in the metric $g$ applying to vector
fields and $\Ric(\cdot,\cdot)$ is the Ricci curvature tensor. Since
$\Ric=\kappa g$ and $\<{\bf \Delta}Dw,Dw\>=-\<D(\Delta w),Dw\>,$ we
have, by (\ref{8.12}) and (\ref{8.15}), \be
|D^2w|_{T_x^2}^2=\frac{1}{2}\Delta|Dw|^2+\kappa|Dw|^2\qfq
x\in\Om.\label{8.16}\ee

From (\ref{8.11}) and (\ref{8.16}), we obtain\beq
|\Xi(V)|_{T_x^2}^2&&=|\kappa
wg+D^2w|_{T_x^2}^2=2\kappa^2w^2+2\kappa w\<g,D^2w\>_{T_x^2}+|D^2w|_{T_x^2}^2\nonumber\\
&&=\frac{1}{2}\Delta|Dw|^2+\kappa(|Dw|^2-2\kappa
w^2)\nonumber\\
&&=\frac{1}{2}\Delta|Dw|^2+\kappa\div w\nabla w \qfq
x\in\Om.\nonumber\eeq     \hfill$\Box$\\

Let $\V_0^\bot(\Om(a))$ be given in (\ref{nn4.11}). By Theorem
\ref{t4.6}, $\V_0^\bot(\Om(a))=H^1_\ib(\Om(a))$ for
$0<a<\dfrac{\pi}{2\sqrt{\kappa}}.$ We define a linear operator
$\Theta:$ $L^2(\Ga(a))\rightarrow L^2(\Ga(a))$ by
$$\Theta\psi=w_\rho,$$ where $w\in \V_0^\bot(\Om)$ is the solution
to the problem (\ref{8.12}). Then $D(\Theta)=H^{1/2}(\Ga(a))$ for
$0<a\leq \dfrac{\pi}{\sqrt{\kappa}}.$

\begin{thm}\label{t8.2} Let $\Om(a)$ and $\Ga(a)$ be given in $(\ref{8.17}).$
Then the bending energy $(\ref{8.3})$ of the $\Ga$-convergence
becomes the following one-dimensional problem \be
\t{I}(\psi)=\frac{\mu}{12}\int_{\Ga(a)}[2\psi_\tau
(\Theta\psi)_{\tau}-\kappa \psi
\Theta\psi-\sqrt{\kappa}a\ctg(\sqrt{\kappa}a)(|\Theta\psi|^2+|\psi_\tau|^2)]d\Ga\label{8.18}\ee
for $\psi\in H^{1/2}(\Ga(a)),$ where $\tau$ is the unit tangential
vector field along $\Ga(a).$
\end{thm}

{\bf Proof}\,\,\,Let $\tau=\tau(\rho)$ be the unit tangential vector
field along $\Ga(\rho)$ for $0<\rho\leq a.$ Then $D\rho,$ $\tau$
forms a frame field on $\Om(a).$  We have \be D_{D\rho}D\rho=0,\quad
D_{D\rho}\tau=0,\label{8.19*}\ee\be D_\tau
D\rho=\sqrt{\kappa}\rho\ctg(\sqrt{\kappa}\rho)\tau,\quad
D_\tau\tau=-\sqrt{\kappa}\rho\ctg(\sqrt{\kappa}\rho)D\rho.\ee
Moreover, the equation in (\ref{8.12}) gives \beq
D^2w(D\rho,D\rho)&&=-2\kappa
w-D^2w(\tau,\tau)\nonumber\\
&&=-w_{\tau\tau}-\sqrt{\kappa}a\ctg(\sqrt{\kappa}a)w_\rho-2\kappa
w\qfq x\in\Ga(a).\label{8.21}\eeq It follows from the formulas
(\ref{8.14}) and (\ref{8.19*})-(\ref{8.21}) that
\beq\int_{\Om(a)}|\Xi(V)|_{T^2_x}^2dg&&=\int_{\Ga(a)}[D^2w(D\rho,Dw)+\kappa
ww_\rho]d\Ga\nonumber\\
&&=\int_{\Ga(a)}[w_\rho
D^2w(D\rho,D\rho)+w_\tau(w_{\rho\tau}-\<Dw,D_\tau D\rho\>)+\kappa
ww_\rho]d\Ga\nonumber\\
&&=\int_{\Ga(a)}[2w_\tau
w_{\rho\tau}-\sqrt{\kappa}a\ctg(\sqrt{\kappa}a)(w_\rho^2+w_\tau^2)-\kappa
ww_\rho]d\Ga.\label{8.22}\eeq

Finally, we use the formulas (\ref{8.22}), (\ref{8.13}) and
(\ref{8.7}) in the formula (\ref{8.3}) to obtain (\ref{8.18}).
\hfill$\Box$\\

{\bf Bending of a Cylinder Shell}\,\,\,Let $a>0$ and let
\be\Om=\{\,\,(\cos\theta,\sin\theta,z)\,\,|\,\,\theta\in[-\pi,\pi),\,\,|z|<a\,\,\}.\label{8.23}\ee
Then
$$\pl z=(0,0,1),\quad\pl\theta=(-\sin\theta,\cos\theta,0).$$

 Let $w\in H^1_\is(\Om)$ be given. By Theorem \ref{t6.1},
$$w=w_0+w_1z,\quad w_0,\,\,w_1\in H^1({\bf T}).$$ Let $W\in\X(\Om)$
be such that $V=W+wN$ is an infinitesimal isometry. A simple
computation shows that
$$   W=[\int_0^\theta(\theta-\eta)w_1(\eta)d\eta+c_1]\pl
z-[\int_0^\theta [w_0(\eta)+w_1(\eta)z]d\eta+c_2]\pl \theta,$$ where
$c_1,$ $c_2$ are constants and
$$\Xi(V)(\pl z,\pl z)=0,\quad \Xi(V)(\pl z,\pl\theta)=-\int_0^\theta w_1(\eta)d\eta-w_{1\theta},$$
$$\Xi(V)(\pl\theta,\pl\theta)=-w-w_{\theta\theta}.$$

Using the above formulas, we obtain

\begin{thm}\label{t8.3} Let $\Om$ be given by $(\ref{8.23}).$ Then the bending energy $(\ref{8.3})$ of the $\Ga$-convergence
becomes the following one-dimensional formula \beq
\t{I}(w_0,w_1)&&=\int_{-\pi}^\pi\Big\{\frac{\mu
a}{3}[\frac{\mu+\lam}{2\mu+\lam}(w_0+w_{0\theta\theta})^2+(w_{1\theta}+\int_0^\theta
w_1(\eta)d\eta)^2]\nonumber\\
&&\quad+\frac{\mu(\mu+\lam)a^3}{3(2\mu+\lam)}(w_1+w_{1\theta\theta})^2\Big\}d\theta\qfq
(w_0,w_1)\in H^1({\bf T})\times H^1({\bf T}).\nonumber\eeq

\end{thm}

\end{document}